\documentclass{amsart}
\usepackage{latexsym}
\usepackage{capt-of}
\usepackage{graphicx}

\DeclareFontFamily{U}{wncy}{}
\DeclareFontShape{U}{wncy}{m}{n}{<->wncyr10}{}
\DeclareSymbolFont{mcy}{U}{wncy}{m}{n}
\DeclareMathSymbol{\Sha}{\mathord}{mcy}{"58}

\usepackage{amssymb,amsmath}
\usepackage{amsfonts}
\usepackage{hyperref}
\newcounter{ctfig}

\newcommand{\Z}{{\mathbb Z}}
\newcommand{\Q}{{\mathbb Q}}

\newcommand{\V}{\mathcal{V}}
\newcommand{\C}{\mathcal{C}}

\DeclareMathOperator{\Kum}{{\rm Kum}}
\newcommand{\dsKum}{\widetilde \Kum}

\renewcommand{\O}{\mathcal{O}}

\newcommand{\A}{\mathbb{A}}

\newcommand{\F}{\mathbb{F}}
\renewcommand{\P}{\mathbb{P}}

\renewcommand{\L}{\mathcal{L}}

\newcommand{\ethreetwo}{{{}_3E_2}}
\newcommand{\fthreetwo}{{{}_3F_2}}
\newcommand{\athreetwo}{{{}_3A_2}}

\DeclareMathOperator{\Spec}{{\rm Spec}}
\DeclareMathOperator{\Aut}{{\rm Aut}}

\DeclareMathOperator{\codim}{{\rm codim}}
\newcommand{\JacC}{{\hbox{Jac}_{\lower.5pt\hbox{$_\C$}}}}
\newcommand{\JacF}{{\hbox{Jac}_{\lower.5pt\hbox{$_\F$}}}}
\newcommand{\etale}{{\textnormal{\' et}}}

\DeclareMathOperator{\tr}{\rm tr}

\DeclareMathOperator{\Sym}{\rm Sym}

\newcommand{\la}{\lambda}
\newcommand{\qr}[2]{\left(\frac{#1}{#2}\right)}

\DeclareMathOperator{\Pic}{\rm Pic}

\newtheorem{definition}{Definition}
\newtheorem{remark}{Remark}
\newtheorem{conjecture}{Conjecture}
\newtheorem{theorem}{Theorem}
\newtheorem{example}{Example}
\newtheorem{proposition}{Proposition}
\newtheorem{corollary}{Corollary}
\newtheorem{lemma}{Lemma}

\bibstyle{plain}
\newcommand\faketheorem[2]{\medskip\noindent{\bf #1. }{\em #2}\medskip}
\begin{document}
%% \foreach \x in{graphics,floats}{%
%%     \immediate\write18{pdflatex -jobname=template-\x\space "\def\noexpand\placeholder{\x} \noexpand\input{template}"}%
%%     \includepdf[pages=-]{template-\x}%
%% }

\title[Modularity of two double covers of $\P^5$]{Modularity of two double covers of $\P^5$ branched along $12$ hyperplanes}
%\titlerunning{Modularity of two double covers of $\P^5$}
%\title{Modularity of two double covers of $\P^5$ branched along $12$ hyperplanes}

\author{Adam Logan}
%% \institute{\author{A. Logan} \at
%%   The Tutte Institute for Mathematics and Computation,
%%   P.O. Box 9703, Terminal, Ottawa, ON K1G 3Z4, Canada;\\
%%   School of Mathematics and Statistics, 4302 Herzberg Laboratories,
%%   1125 Colonel By Drive, Carleton University, Ottawa, ON K1S 5B6, Canada\\
%%   \email{adam.m.logan@gmail.com}}
\address{The Tutte Institute for Mathematics and Computation,
P.O. Box 9703, Terminal, Ottawa, ON K1G 3Z4, Canada}
\address{School of Mathematics and Statistics, 4302 Herzberg Laboratories,
  1125 Colonel By Drive, Carleton University, Ottawa, ON K1S 5B6, Canada}

\begin{abstract}
  For two varieties of dimension $5$ constructed as double covers of $\P^5$
  branched along the union of $12$ hyperplanes, we prove that the number of
  points over~$\F_p$ can be expressed in terms of Artin symbols and the
  $p$th Fourier coefficients of modular forms.  Many analogous results are
  known in dimension $\le 3$, but very few in higher dimension.  In addition,
  we use an idea of Burek to construct quotients of our varieties for which the
  point counts mod $p$ are expressible in terms of Artin symbols
  and the coefficients of a single modular form of weight~$6$.
\end{abstract}
\date{\today}

\maketitle

\section{Introduction}\label{sec:intro}
One of the most fundamental problems in number theory is to give formulas for
the number of points on a fixed variety over a varying finite field.
In the simplest cases,
the number of points over~$\F_p$ can be expressed in terms of powers of $p$
and Artin symbols expressing the decomposition of $p$ in number fields;
however,
this is very far from being sufficient in general.  Most famously, isogeny
classes of elliptic curves over $\Q$ of conductor $N$
are in bijection with newforms of level~$N$ with integer coefficients
\cite{bcdt};
the correspondence takes a form with Hecke eigenvalues $a_p$ to an elliptic
curve with $p+1-a_p$ points over~$\F_p$.

In higher dimensions, we cannot expect such a statement to hold for all
varieties of a given deformation type, except in very special situations.
Nevertheless, we would like to find varieties for which the number of points
can be expressed in terms of powers of $p$, Artin symbols, and the
coefficients of
modular forms.  Perhaps the most interesting case is that in which only a
single eigenform of weight~$n>2$ is needed and the dimension of the variety
is $n-1$.

In dimension greater than $1$, the most natural candidates for this property
are the {\em Calabi-Yau varieties}, which are defined as follows:

\begin{definition}\cite[Definition 1]{gouvea-yui}\label{defn:calabi-yau}
  A {\em Calabi-Yau variety} is a smooth 
  variety $V$ of dimension $d$ satisfying
  $K_V \cong \O_V$ and $H^i(K_V) = 0$ for $0 < i < d$.  If $V$ is a 
  limit of Calabi-Yau varieties and has a Calabi-Yau resolution of
  singularities, then $V$ is a {\em singular Calabi-Yau variety}.
  If $H^{d-1}({\mathcal T_V}) = 0$, where ${\mathcal T_V}$ denotes the tangent
  bundle of $V$, then $V$ is {\em rigid}; this condition can also be written
  as $H^{d-1,1}(V) = 0$.  If $h^{i,j} = 0$ for all $i,j$ except with $i = j$
  or $\{i,j\} = \{0, \dim V\}$, and either $\dim V$ is odd or the
  semisimplification of $H_\etale^{\dim V}(V,\Z_p)$ splits
  as a direct sum of Galois representations
  $(H^{\dim V,0} \oplus H^{0,\dim V}) \oplus H^{\dim V/2,\dim V/2}$,
  then $V$ is {\em strongly rigid}.
\end{definition}

\begin{remark}\label{rem:basic-cys}
  According to our definition, an elliptic curve is a strongly rigid
  Calabi-Yau variety (sometimes these are defined to require trivial
  fundamental group).
  A K3 surface is a Calabi-Yau variety, and it is strongly rigid if and only
  if its Picard number is $20$.  A Calabi-Yau threefold is rigid if and only
  if it is strongly rigid.
  
  One can prove that a Calabi-Yau variety $V$ is strongly rigid by
  showing that the cohomology of $V$ is generated by cycle classes and
  $H^{\dim V,0}, H^{0,\dim V}$.  If $\dim V$ is odd this is equivalent,
  and presumably that is true if $\dim V$ is even, but a proof would seem to
  need some form of the Hodge conjecture.
\end{remark}
  
The integral Hecke eigenforms of weight~$3$ up to twist correspond to imaginary
quadratic fields with class group of exponent dividing $2$ by
\cite[Theorem 2.4]{schutt}.
Such fields are known with at most one exception, which is excluded by
the generalized Riemann hypothesis; the list can be found 
in \cite{elkies-schutt}.  Elkies and Sch\" utt proved the following:

\begin{theorem}[{\cite[Theorem 1]{elkies-schutt}}]\label{thm:e-s}
  Let $f$ be a Hecke eigenform of weight~$3$ from the list
  with eigenvalues $a_p$.  Then there is a K3 surface $S_f$
  such that, for all but finitely many $p$, we have
  $\#S_f(\F_p) = p^2 + c(p)p + 1 + a_p$, where $c(p)$ is a linear combination of
  Artin symbols.
\end{theorem}

Much less is known for forms of weight greater than $3$.
Gouv\^ea and Yui proved the following:
\begin{theorem}[{\cite[Theorem 3]{gouvea-yui}}] Let $V$ be a smooth rigid
  Calabi-Yau threefold over $\Q$.  Then $V$ is modular: in other words,
  there is a Hecke eigenform $f$ of weight~$4$ such that the Galois
  representation $\rho_\ell$ on $H^3_\etale(V,\Q_\ell)$ is equivalent to
  $\rho_{f,\ell}$, the representation attached to $f$, for all $\ell$.
\end{theorem}
(The same conclusion had been obtained by
Dieulefait and Manoharmayum \cite{dieu-man} under stronger hypotheses.)
It follows that $V$ has $p^3 + n(p)(p^2+p) + 1 - a_p$ points over~$\F_p$
for all primes $p$ of good reduction,
where $n(p)$ is expressed in terms of Artin symbols and the $a_p$ are
the Hecke eigenvalues of $f$.
Many examples are worked out in detail in \cite{meyer} and elsewhere,
but it is not known
whether every Hecke eigenform can be realized by a rigid Calabi-Yau threefold
in this way, nor whether there are finitely or infinitely many rational
Hecke eigenforms of weight~$4$ up to twist.

Gouv\^ea and Yui point out that their methods also apply to Calabi-Yau
varieties of odd dimension $d$ with $\dim H^d = 2$
\cite[page~146]{gouvea-yui}.
Thus, if $V$ is a strongly rigid Calabi-Yau
variety of odd dimension $d$, there is a formula for the number of
$\F_p$-points
of $V$ of the form $\sum_{i=1}^d n_i p^i + (-1)^d a_p$, where the $n_i$ are
expressed in terms of Artin symbols and the $a_p$ are the Hecke eigenvalues
of a rational newform of weight~$d+1$, and the same is expected to hold
for $d$ even.  On the other hand, if $V$ is only
rigid, this would not be expected.  For example, if $d = 5$ it is possible
that $h^{3,2}(V) > 0$ and that
$H^5_{\etale}(V,\Z_\ell)$ is irreducible of dimension $>2$ even though $V$ is
rigid.

In dimension greater than $3$ there are almost no examples of strongly
rigid Calabi-Yau varieties.  If the dimension of the
space of cusp forms of weight~$k$ for $\Gamma_1(N)$ is $1$, then the Kuga-Sato
construction \cite{deligne} gives a Calabi-Yau variety realizing the form.
Ahlgren \cite{ahlgren} in effect studies the case $k = 6, N = 4$,
and Paranjape and Ramakrishnan \cite{p-r} consider several others, including
$k = 6, N = 3$.
In addition, Frechette, Ono, and Papanikolas have
shown \cite{fop} how to construct varieties that realize the cusp forms of
level $N = 2, 4, 8$ and arbitrary weight~$k$.  However, for $k$ even and
greater than $10, 6, 4$
respectively, the dimension of the space of cusp forms of level~$N$ and weight
$k$ is greater than $1$ and the variety they construct is not Calabi-Yau.
Roberts has conjectured \cite[Conjecture 1.1]{roberts} that up to twist there
are only finitely many newforms with rational coefficients and not of CM type
(complex multiplication, i.e., for which there is a nonsquare integer $N$
such that
$a_p = 0$ for all primes $p$ with $\qr{N}{p} = -1$) for $k \ge 6$ and none for
$k \ge 52$.

The main goal of this paper is to work out two examples of fivefolds.  One
realizes the newform of weight~$6$ and level~$8$; the other, the newform
of weight~$6$ an level~$32$ with complex multiplication.
Both are double covers
of $\P^5$ branched along a union of $12$ hyperplanes; however, the methods are
somewhat different.

Our analysis of the first example, which bears a close resemblance to certain
rigid Calabi-Yau threefolds, will culminate in the proof of the following
result:

\faketheorem{Theorem~\ref{thm:main-first}}{The
  double cover $V_8$ of $\P^5$ defined by
  $t^2 = \prod_{i=0}^5 x_i (x_i+x_{i+1})$ is modular
  (Definition~\ref{def:modular}) of level~$8$.}

We do this by 
finding a fibration by quotients of products of Kummer surfaces closely
related to the construction of~\cite{fop}.  This variety appears not to be
a singular Calabi-Yau, but the modularity does
not require such a statement to hold in any case.
We will also use an idea of Burek \cite{burek}
to construct a quotient of the variety that appears to be a strongly rigid
Calabi-Yau.  That is, we will prove:

\faketheorem{Theorem~\ref{thm:q3}}{Let $Q_3 = V_8/\iota_3$,
    where $\iota_3$ takes $x_i$ to $x_{6-i}$.
    Then the number of $\F_p$-points of $Q_3$ is
    $\sum_{i=0}^5 p^i - a_p - \phi(-1)p^2$ for all odd $p$,
    where $a_p$ is the $p$\/th Fourier coefficient of the newform of
    weight~$6$ and level~$8$ and $\phi$ is the quadratic character mod $p$.}

The second example is also a double cover of $\P^5$ branched along
the union of $12$ hyperplanes.  We summarize our results on it in the
following statements:

\faketheorem{Theorem~\ref{thm:count-32}, Proposition~\ref{prop:bir-v32},
  Theorem~\ref{thm:rigid-32}}{Let $V_{32}$ be the fivefold defined by
  $$t^2 = \left(\prod_{i=0}^5 x_i\right)(x_0+x_1)(x_3+x_5)(x_2+x_4+x_5)(x_0+x_2-x_4)(x_1-x_2+x_4)(x_2-x_3+x_4).$$
  Let $E$ be the elliptic curve $y^2 = x^3 - x$ and $M$ the K3 surface
  defined by $t^2 = xyz(x+y)(y+z)(-x+z)$.  Then $V_{32}$ is birationally
  equivalent to the quotient of $M \times M \times E$ by a group of
  order $4$.  Further, for all odd primes $p$, there are
  $\sum_{i-0}^5 p^i - a_{6,p} - pa_{4,p} - 2p^2a_{2,p}$ points on $V_{32}$ over
  $\F_p$, where $a_{i,p}$ is the $p$th Fourier coefficient of the newform
  of weight~$i$ and level~$32$ that has complex multiplication.  A
  quotient of $V_{32}$ by a group of order $4$ is birationally equivalent
  to a variety with $\sum_{i=0}^5 p^i - a_{6,p}$ points over~$\F_p$ for all
  odd $p$.}

Here $M$ has Picard
number $20$ and realizes the CM newform of weight~$3$ and level~$16$ with
quadratic character.  We will also indicate
a related construction in level~27 in Example~\ref{ex:level-27}.
%and in Example~\ref{ex:several-forms} we will show how to generalize it
%to involve modular forms of different levels.

An extensive search has discovered many examples of apparently modular
double covers of $\P^5$ with branch locus the union of $12$ hyperplanes
that will not be discussed in this paper.
One of these is Ahlgren's fivefold \cite{ahlgren}; others
appear to be new and we intend to study them in future work.
In particular, we have found two more unions of
$12$ hyperplanes, not projectively equivalent to the
first one or to each other, such that the double cover has a Calabi-Yau
resolution that appears to be strongly rigid and to realize the form of
level~$8$.  One of them is notable for its large symmetry group, with the
symmetric group on $5$ symbols acting faithfully on the set of $12$
hyperplanes; the other, for admitting a fibration in quotients of products of
K3 surfaces that is similar to
but distinctly different from that in the first example discussed here.

The conjectural identity that equates the number of points on this variety
computed from the fibration to a simple formula in terms of the coefficients
of the modular form appears to point to a previously undiscovered identity of
hypergeometric functions.  Similarly, by considering quotients of our
example of level~$32$ and comparing the point counts that arise from
fibrations to those obtained in this paper and from its results, one obtains
some attractive identities, which again will be presented in future work.
There are also several collections of hyperplanes that appear to
correspond to a cusp form of level~$256$ which is a quartic twist of the
form of level~$32$.  It is interesting that we do not find any examples whose
level is not a power of $2$ (in dimension $3$ there are many such examples
of double covers of $\P^3$ branched along the union of $8$~hyperplanes).

\subsection{Acknowledgements}
%\begin{acknowledgements}
  I thank the referees for their helpful comments, suggesting
  among other things a substantial reorganization of the contents of the
  paper into a more logical order.  I also thank Colin Ingalls, Ken Ono,
  Owen Patashnick, Rachel Pries, John Voight, and Don Zagier for enlightening
  discussions, and Spencer Secord for his remarks on a draft of this paper.
%\end{acknowledgements}

\section{Notation}\label{sec:notation}
We start by introducing some notation that will apply throughout the paper.

\begin{definition}\label{def:notation}
  We will often work in $\P^5_\Q$; the 
  coordinates in it will be denoted by $x_0, \dots, x_5$,
  with the usual understanding that $x_i = x_{i+6}$.  We will also use
  weighted projective space with weights $6,1,1,1,1,1,1$, this being the
  natural home for double covers of $\P^5$ with branch locus of degree $12$.
  The variables there will be $t, x_0, \dots, x_5$, and a map from
  $\P(6,1,1,1,1,1,1)$ to $\P^5$ will always be given by omitting $t$.
  At times we will use
  other projective spaces, referring to their variables as $y_0, \dots, y_m$
  or $z_0, \dots, z_n$.  The weighted projective space $\P(3,1,1,1)$ will
  also arise, since our arguments use many K3 surfaces given as double covers
  of $\P^2$; its coordinates will usually be $t, x, y, z$, but sometimes
  $w, z_0, z_1, z_2$, etc.
\end{definition}

\begin{definition} Let $B$ be a scheme with a given $\F_p$-point
  $\iota: \Spec{\F_p} \hookrightarrow B$ and
  $\pi: \V \to B$ a flat family of schemes of finite type over $B$.
  Then $\Spec{\F_p} \times_B \V$ is a scheme $V$ of finite type over
  $\F_p$, so it has finitely many $\F_p$-rational points.  We denote
  the number of these points by $[\V]_p$ (the choice of $\iota$ will always
  be clear).
  By abuse of language we will also use this notation when $\V$
  is a variety defined over $\Q$ by equations
  with coefficients whose denominators are not multiples of $p$,
  implicitly viewing it as a scheme over $\Spec \Z_{(p)}$.
\end{definition}

\begin{definition}\label{def:modular}
  Let $V$ be a variety over $\Q$.  If there is a formula
  for $[V]_p$, valid for all but finitely many $p$, in terms of powers of
  $p$, Artin symbols, and eigenvalues of Hecke eigenforms for $\Gamma_0(N)$,
  then $V$ is {\em modular}, and its {\em level} is the least common multiple
  of those of the eigenforms that appear.
  In the opposite direction we say that an
  eigenform involved in such a formula is {\em realized} by $V$,
  especially if it is the only eigenform of its weight that appears in the
  formula.
  %% We could of course consider a broader concept of modularity
  %% that permits Hilbert modular forms, Siegel modular forms, etc., but
  %% such objects will not appear here.
\end{definition}

\begin{definition} We use $\phi$ to denote the quadratic character modulo a prime
  $p$ (we will never be considering more than one prime at a time, so this
  will not lead to ambiguity).
\end{definition}

\section{Counting points on double covers}
In this paper we will have to count the $\F_p$-points of various double
covers of affine and projective spaces and of quotients of products of
these.  We give the notation and state the results that we will be using
here so as not to
interrupt the exposition later.  The main result of this section,
Lemma~\ref{lem:count-q-dc}, is certainly well known, but we prove it for
lack of an appropriate reference.

Let $S$ be affine or projective space over~$\F_p$ and let $f$ be a
polynomial defining a subvariety of $S$; if $S$ is projective we suppose
further that $\deg f$ is even.  Then there is a double cover $D_f$ of $S$
defined by the equation $s^2 - f = 0$, where $s$ is a new variable.  It is
trivial to compute the number of $\F_p$-points of $D_f$ in terms of the
values of $f$:

\begin{lemma}\label{lem:count-one-dc}
  We have $[D_f]_p = \sum_{P \in S(\F_p)} \phi(f(P))$.
\end{lemma}

(If $S$ is projective then $f(P)$ is not well-defined, but $\phi(f(P))$
still is.)  More generally, we can count points on quotients of products
of double covers by the involution that is the product of the involution
on each factor.  To do so, we first introduce some notation:
\begin{definition}\label{def:count-double-cover}
  Let $V$ be a variety over a finite field $F$ with an involution $\iota$.
  For a positive integer $i$ let $F_i$ be an extension of $F$ of degree $i$.
  Let $P_{V,i}$ (respectively~$Z_{V,i}, N_{V,i}$) be the number of $F_i$-points of
  $V/\iota$ whose inverse images in $V/F_i$ contain $2$ (resp.~$1, 0$)
  rational points.  This depends on $\iota$, but there will never be any
  ambiguity, so we suppress it from the notation.  We refer to a point $P_0$
  of $V/\iota$ as a {\em $P$-point} (resp.~{\em $Z$-point},
  {\em $N$-point}), depending on the number of points of $V$ above $P_0$.
\end{definition}

\begin{lemma}\label{lem:count-q-dc}
  Let $p$ be an odd prime.  For $1 \le i \le n$, let $S_i$ be affine or
  projective spaces over~$\F_p$ and let $V_i$ be hypersurfaces defined by
  $f_i = 0$ in $S_i$, where the degree of $f_i$ is even if $S_i$ is
  projective.  
  Let the $D_{f_i}$ be the double covers of $S_i$ defined by $s_i^2 - f_i = 0$.
  Then $(\Z/2\Z)^n$ acts on $\prod_{i=1}^n D_{f_i}$ by negating appropriate
  subsets of the $s_i$; let $E$ be the subgroup
  of $(\Z/2\Z)^n$ consisting of elements of even weight,
  and let $D_V = (D_{f_1} \times \dots \times D_{f_n})/E$.  Then
  $[D_V]_p = \prod_{i=1}^n [S_i]_p  + \prod_{i=1}^n(P_{V_i,1} - N_{V_i,1})$.
\end{lemma}

\begin{proof}
  Let $R = (R_1,\dots,R_n)$ be a point of $S_{V_1} \times \dots S_{V_n}$.
  First, if at least one of the $R_i$ is a $Z$-point, then
  $\prod_{i=1}^n f_i(R_i) = 0$ so $(R_1,\dots,R_n)$ lies under one
  rational point of $D_V$.  Otherwise there are $2^n$ geometric points
  of $\prod_{i=1}^n D_{f_i}$ lying over $R$, and a point and its Galois
  conjugate are in the same $E$-orbit if and only if an even number of
  the $R_i$ are $N$-points.  In this case there are $2$ rational points
  above $R$ in $D_V$, and otherwise there are $0$.  Equivalently,
  the number of $\F_p$-points of $D_V$ above $R$ is
  $1 + \prod_{i=1}^n (p_i(R_i)) - n_i(R_i))$, where $p_{P_i}, n_{P_i}$ are the
  characteristic functions of the sets of $P$- and $N$-points on $V_i$.
  The result follows by summing over all $\F_p$-points of
  $S_1 \times \dots \times S_n$.
\end{proof}

\section{Hypergeometric functions over finite fields and modular forms}\label{sec:hypergeom}
The main purpose of this brief section is to relate our notation for elliptic
curves to that used by Frechette, Ono, and Papanikolas in
\cite{fop} so as to apply their results
relating hypergeometric functions over finite fields to Hecke eigenvalues
of modular forms.  We begin by recalling their notation.

\begin{definition}[{\cite[(2.2), (2.4), (1.1), (1.3)]{fop}}]\label{fop-notation}
  Let $\ethreetwo(\la)$ be the elliptic curve defined by
  $y^2 = (x-1)(x^2+\la)$ and, for $\la \in \F_p$ with
  $\la^2 \ne -\la$, let $\athreetwo(p,\la)$ be the trace of
  Frobenius of $\ethreetwo(\la)$ over~$\F_p$.  In addition, for characters
  $A$ and $B$ on $\F_p$, let $\binom{A}{B}$ be the normalized Jacobi sum
  $\frac{1}{p}\sum_{x \in \F_p} A(x) \bar B(x-1)$, where the bar denotes the
  complex conjugate.
  Let $\phi$ be the quadratic character on $\F_p$ and let
  $\fthreetwo(\la) = \frac{p}{p-1} \sum_\chi {\binom{\phi \chi}{\chi}}^3 \chi(\la)$,
  where the sum runs over all characters $\chi$ of $\F_p$.
\end{definition}

We restate the basic relation between $\fthreetwo$ and $\athreetwo$.

\begin{theorem}[{\cite[Theorem 4.3 (2), Theorem 4.4 (2)]{fop}}]
  $$\fthreetwo\left(1+\frac{1}{\la}\right) = \frac{\phi(-\la)(\athreetwo(p,\la)^2-p)}{p^2}$$
  for $\la \in \F_p$ with $\la \ne 0, -1$.  In addition, if
  $p \equiv 1 \bmod 4$ we have $\fthreetwo(1) = \frac{4a^2-2p}{p^2}$ where $a$ is
  an odd integer such that $p - a^2$ is a square, and if $p \equiv 3 \bmod 4$
  then $\fthreetwo(1) = 0$.
\end{theorem}

To relate our notation to that of~\cite{fop} requires a simple statement
about elliptic curves.

\begin{definition}\label{def:e-lambda} For $\la \ne 0, -1$, let
  $E_\la$ be the elliptic curve defined by
  $y^2 = x^3 -2x^2 + \frac{\la}{\la+1}x$.  
  Let $a_{\la,p} = p+1-\#E_\la(\F_p)$, and let $a_{0,p} = 0$.
\end{definition}

\begin{proposition}\label{prop:same-curve} For $\la \ne 0, -1$ we have
  $\athreetwo(p,\frac{-1}{\la+1})^2 = a_{\la,p}^2$.  Equivalently, we have
  $\athreetwo(p,\mu)^2 = a_{-(1+\frac{1}{\mu}),p}^2$ for $\mu \ne 0, -1$.
\end{proposition}

\begin{proof} Replacing $x$ by $x+1$ in the equation
  $y^2 = (x-1)(x^2-\frac{1}{\la+1})$ defining an elliptic curve whose trace
  of Frobenius is $\athreetwo(p,\frac{-1}{\la+1})$ gives a quadratic
  twist of the elliptic
  curve $y^2 = x^3-2x^2+\frac{\la}{\la+1}x$.  This is the elliptic curve whose
  trace is $a_{\la,p}$, so the two have the same trace up to sign.
\end{proof} 

Combining these two statements gives
\begin{equation}\label{f-a}
  \fthreetwo(\la) = \phi\left(\frac{\la+1}{\la}\right)\frac{(a_{-\la,p}^2-p)}{p^2}.
\end{equation}

The following simple statements about modular forms with
complex multiplication by $\Q(i)$ will be used in
the proof of Theorem~\ref{thm:count-32}.

\begin{definition}\label{def:mj-ajp}
  For $j \in \{2,4,6\}$, let $m_j$ be the unique newform of weight
  $j$ and level~$32$ that has complex multiplication by $\Q(i)$.  Let
  $m_3$ be the newform of weight~$3$ and level~$16$ whose Nebentypus is the
  Dirichlet character $\left(\frac{-1}{\cdot}\right)$.  (In the LMFDB
  \cite{lmfdb} these are {\tt 32.2.a.a}, {\tt 32.4.a.b}, {\tt 32.6.a.b},
  and {\tt 16.3.c.a} respectively.)  For $j \in \{2,3,4,6\}$
  and $p$ prime, let $a_{j,p}$ be the eigenvalue of $m_j$ for the Hecke
  operator $T_p$.
\end{definition}

The following is a routine application of the theory of modular forms with
complex multiplication.
\begin{lemma}\label{lem:coef-mf} For $p \equiv 1 \bmod 4$ we have
  $a_{3,p} = a_{2,p}^2 - 2p,
  a_{4,p} = a_{2,p}(a_{3,p}-p), 
  a_{6,p} = a_{4,p}a_{3,p} - p^2 a_{2,p}$; for $p \equiv 2, 3 \bmod 4$ all
  $a_{j,p}$ are equal to $0$.
\end{lemma}

\begin{proof}
  Since the $m_i$ are modular forms with complex multiplication by $\Q(i)$,
  their Fourier coefficients may be described in terms of Hecke characters
  of this field \cite{ribet}, and in particular $a_{j,p} = 0$ if $p$ does not
  split in $\Q(i)$.
  For $p \equiv 1 \bmod 4$, let $a_p, b_p$
  be such that $a_p^2 + b_p^2 = p$ and $a_p + b_p \equiv 1 \bmod 2+2i$
  (this determines $a_p$ uniquely and $b_p$ up to sign).  Then
  $a_{j,p} = \tr (a_p + b_p i)^{j-1}$ for $j \in \{2,3,4,6\}$.  The claim is now
  easily checked.
\end{proof}

%% \begin{proof}
%%   For the second identity
%%   $$\begin{aligned}
%%     a_{4,p} &= \tr(a_p+b_p i)^3 \cr
%%     &= 2a_p^3 - 6a_pb_p^2\cr
%%     &= a_{2,p}(a_p^2 - 3b_p^2)\cr
%%     &= a_{2,p}(2a_p^2 - 2b_p^2 - p).\cr
%%   \end{aligned}$$
%%   But $a_{3,p} = \tr(a_p+b_pi)^2 = 2a_p^2 - 2b_p^2$, so the first claim
%%   follows.  The proofs of the other two are similar.
%%   % &a_{6,p} = \tr(a_p+b_p i)^5 \cr
%%   % &= 2a_p^5 - 20a_p^3 b_p^2 + 10a_p b_p^4\cr
%%   % &= a_{2,p}(a_p^4 - 10a_p^2b_p^2 + 5b_p^4)\cr
%%   % so we need to show that this is (a_{3,p}-p)(a_{3,p}) - p^2.
%%   % well, a_{3,p} = 2a_p^2 - 2b_p^2 so a_{3,p}-p = a_p^2 - 3b_p^2
%%   % and the product is 2a_p^4 - 8a_p^2 b_p^2 + 6bp^4.
%%   % Subtract (a_p^2 + b_p^2)^2 to get what we want.
%% \end{proof}

\section{Geometry and point counting on K3 surfaces}\label{sec:geom-points}
In this section we prove some results about the Kummer surface
$(E \times E)/\pm 1$ of the
square of an elliptic curve and some closely related elliptic surfaces.
These will be used repeatedly in the remainder
of the paper.  While some of these statements are not easily found in the
literature, they are all fairly routine exercises, so few detailed proofs
will be given.  We begin by introducing two surfaces that are isogenous to
the Kummer surface and a quadratic twist, which are important in
Section~\ref{sec:first}.

\begin{definition}\label{def:kla-lla}
    Let $K_\la, L_\la$ be the surfaces in 
  $\P(3,1,1,1)$ defined by the equations
  \begin{equation}
    \label{kla-lla}
    \begin{split}
      v^2 &= (\la+1) z_0z_1z_2(\la z_0+z_1)(z_1+z_2)(z_0+z_2), \\
      w^2 &= \la(\la+1)y_0y_1y_2(y_0+y_1)(\la y_0+y_2)(y_1+y_2)\\
    \end{split}
  \end{equation}
  and let
  $A_\la, B_\la$ be the affine patches $z_0 \ne 0, y_0 \ne 0$.
  (The twist by $\la(\la+1)$
  is made to facilitate the comparison with a Kummer surface.)
\end{definition}

\begin{lemma}\label{lem:ka-lb}
  For all $p$ and all $\la \ne 0, -1 \in \F_p$ we have
   $[K_\la]_p - [A_\la]_p = [L_\la]_p - [B_\la]_p = p+1$.
\end{lemma}

\begin{proof}
  The points of $K_\la \setminus A_\la$
  and $L_\la \setminus B_\la$ are exactly the projective points with
  $z_0 = v = 0$.
\end{proof}

By exchanging $z_1, z_2$ we see that $K_\la$ and $L_\la$
are quadratic twists of each other by $\la$.  In other words,
we have $$[K_\la]_p - (p^2+p+1) = \phi(\la)([L_\la]_p - (p^2+p+1)).$$

\begin{definition}\label{def:kum-lambda}
Let $\Kum_\la$ be the Kummer surface
  of $E_\la \times E_\la$: it has a singular model in weighted projective
  space $\P(3,1,1,1)$ defined by
  $$v^2 = \prod_{i=0}^1 (z_i^3-2z_i^2z_2 + \frac{\la}{\la+1}z_iz_2^2).$$
  Let $\dsKum_\la$ be its minimal desingularization.
  
\end{definition}

\begin{proposition}\label{prop:count-kum} Let $\la \ne 0, -1 \in \F_p$.  Then
  $\dsKum_\la$ has
  $p^2+\left(12+6\phi(\la)\right)p+1+a_{\la,p}^2$ points
  over~$\F_p$ (for $a_{\la,p}$ see Definition~\ref{def:e-lambda}).
\end{proposition}

\begin{proof} Consider the elliptic fibration $(z_0:z_2)$ of
  constant $j$-invariant.  If $\la \in \F_p^2$, the $I_0^*$ fibres and their
  components are all rational and contain $4(5p+1)$ points in total.
  In addition, there are $(p-3-a_{\la,p})/2$ fibres isomorphic to $E_\la$
  and containing $p+1-a_{\la,p}$ points and $p-3+a_{\la,p}$ fibres isomorphic
  to the quadratic twist and containing $p+1+a_{\la,p}$ points.  The total
  number of points is $p^2+18p+1+a_{\la,p}^2$.  The argument when
  $\la \notin \F_p^2$ is similar.
\end{proof}

\begin{proposition}\label{prop:fib-on-kum}
  There is a genus-$1$ fibration on $\dsKum_\la$ 
  the Jacobian of whose general fibre is isomorphic to the
  elliptic curve defined by
  $$y^2 = x^3 + \frac{4t-2}{(\la+1)t(\la t+1)}x^2 +
  \frac{1}{((\la+1)(\la t+1)t)^2} x$$
  and that has singular fibres of types $I_4^*, I_1^*, I_0^*, I_1$.
\end{proposition}

\begin{proof}
  We define the fibration by the equations
  $$[z_0 z_1/\la - z_2^2/(\la+1):z_0^2 - 2z_0z_2 + \la z_2^2/(\la+1)].$$
  In Magma \cite{magma}
  it is routine to define the general fibre of this map, verify that
  it is a curve of geometric genus $1$, calculate its Jacobian,
  and show that it has the desired properties.
\end{proof}

\begin{proposition}\label{prop:fib-on-k} The minimal desingularization $\dsKum_\la$
  of $K_\la$ admits a fibration in curves of genus $1$
  whose general fibre is $2$-isogenous to that of the fibration on
  $\Kum_\la$ introduced in Proposition~\ref{prop:fib-on-kum}.  The singular
  fibres of this fibration are of types
  $I_2^*, I_2^*, I_0^*, I_2$, and all components of the
  reducible fibres are defined over the field to which $\la$ belongs.
\end{proposition}

\begin{proof} The desired fibration is defined by $(z_0:z_1)$.
  Again, it is a simple matter to show that the general fibre is isomorphic
  to the elliptic curve defined by
  $$y^2 = x^3 + \frac{t-2}{t(\la+1)(\la t+1)}x^2 +
  \frac{1-t}{((\la+1)(\la t+1)t)^2}x,$$ that its bad fibres are as stated,
  and that the quotient map by the subgroup of order $2$ generated by
  $(\frac{1-t}{(\la+1)(\la t^2+t)}:0:1)$ is the desired isogeny.
\end{proof}

%% \begin{remark}\label{rem:six-lines-isog-kummer}
%%   Note that, up to twist, every elliptic curve over a field of characteristic
%%   not equal to $2$ with a point
%%   of order $2$ and $j$-invariant not equal to $1728$ is a specialization
%%   of $E_\la$.  Indeed, we may write it as $y^2 = x^3 - 2ax^2 + bx$ and twist
%%   by $1/a$, then note that $\la/(\la+1)$ takes all values other than $1$
%%   as $\la$ varies.  Thus we have shown that if $E$ is an elliptic curve with a
%%   point of order $2$, then the Kummer surface of $E \times E$ is isogenous
%%   to a double cover of $\P^2$ branched along six rational lines (this
%%   certainly holds when $j(E) = 1728$ for which this argument does not apply).
%%   This is not really a new result, as it is implicit in~\cite{aop}.
%% \end{remark}

\begin{theorem}\label{thm:count-kl} Suppose as before that $\la \ne 0, -1$.
  Then $[K_\la]_p = p^2 + 1 + a_{\la,p}^2$.
\end{theorem}

\begin{proof}
  We compare the numbers of points on the desingularizations
  by means of the fibrations of Propositions~\ref{prop:fib-on-kum},~\ref{prop:fib-on-k}.
  The number of points on a smooth curve of genus $1$ is unchanged by an
  isogeny, so the difference is accounted for by the singular fibres.
  (The existence of a section is not an issue since every smooth curve of
  genus $1$ over a finite field has rational points.)
  Let $\delta$ be $-1$ if $\Kum_\la$ has split multiplicative reduction
  at the $I_1$ fibre and $1$ otherwise.  
  Then this fibre has $p+1+\delta$ points.
  One checks that the $I_2$ fibre is split if and only if the $I_1$ fibre
  is, so it has $2p+1+\delta$ points.

  All $21$ components of the reducible fibres on $K_\la$ are defined over the
  ground field.  When $\la \in \F_p^2$ the same is true for the $20$ components
  of reducible fibres of $\dsKum_\la$, but
  otherwise only $8$ of them are.
  %: the double curve and two of the tails in the
  %$I_0^*$, the two double curves and two of the tails in the $I_1^*$, and
  %the central component of the $I_4^*$.
  Thus, when $\la \in \F_p^2$, there are $21p+4+\delta$ points on singular
  fibres of the fibration on $\dsKum_\la$
  and likewise $21p+4+\delta$ points on
  singular fibres of the fibration on $\widetilde K_\la$.
  Hence $\widetilde K_\la$ has the same number of points as
  $\dsKum_\la$.

  Similarly, when $\la \notin \F_p^2$, there are $9p+4+\delta$ points on
  singular fibres on $\dsKum_\la$, so $\widetilde K_\la$ has $12p$
  more points than $\dsKum_\la$.  We conclude, in view of Proposition
 ~\ref{prop:count-kum},
  that $\widetilde K_\la$ has $p^2+18p+1+a_{\la,p}^2$ points.

  To finish the proof, we compute that for generic $\la$, the singular
  subscheme of $K_\la$ has degree $18$, and that all components of the
  resolutions are defined over the base field.
  The singular subscheme is unaltered
  by specializations that do not cause additional pairs of lines in the
  ramification locus to meet: the only $\la$ for which such coincidences
  occur are $0, -1$, which are disallowed in the statement of the theorem.
  Hence $\widetilde K_\la$ has
  $18p$ more points than $K_\la$, and the result follows.
\end{proof}

\begin{corollary}\label{cor:count-ll} Let $\la \ne 0, -1$ as before.  Then
  $[L_\la]_p = p^2 + p + 1 + \phi(\la) (a_{\la,p}^2-p)$.
\end{corollary}

\begin{proof} This follows immediately from the fact that $L_\la$ is the twist
  of the double cover $K_\la \to \P^2$ by $\la$.  If $\la$ is a square the two
  surfaces have the same number of points; if not, the sum is twice the number
  of points of $\P^2$.
\end{proof}

We sharpen these results by expressing them in terms of $P$-, $Z$-, and
$N$-points (Definition~\ref{def:count-double-cover}).
This is checked in {\tt count-quotient.mag} \cite{magma-scripts}
by comparing our formulas
to point counts computed directly from equations.

\begin{proposition}\label{prop:pzn-a-b}
  Let $\lambda \ne 0, -1$.  Then over a finite field $F$ of $q$ elements,
  the quantities $P_{A_\la,1}, Z_{A_\la,1}, N_{A_\la,1}$ are equal
  to $(q^2-6q+7+a_{\la,q}^2)/2, 5q-7, (q^2-4q+7-a_{\la,q}^2)$ respectively.
  If $\la$ is a square in $F$, then $P_{B_\la,1}, Z_{B_\la,1}, N_{B_\la,1}$ are
  equal to $P_{A_\la,1}, Z_{A_\la,1}, N_{A_\la,1}$; otherwise they are equal to
  $N_{A_\la,1}, Z_{A_\la,1}, P_{A_\la,1}$.
\end{proposition}

\begin{proof}
  To count $Z_{A_\la,1}$, we note that the branch locus consists of $5$ lines,
  of which $5$ pairs and one triple intersect in rational points.  Thus
  the total number of points is $5q - 5\cdot 1 - 1\cdot 2 = 5q-7$.
  Now $[A_{\la,1}]_q = q^2-q+a_{\la,q}^2$, by Theorem~\ref{thm:count-kl}.
  The result for $A_{\la,1}$ follows by solving the equations
  \begin{align*}
    P_{A_\la,1} + N_{A_\la,1} + Z_{A_\la,1} &= q^2,\\
    2P_{A_\la,1} + N_{A_\la,1}\phantom{+ Z_{A_\la,1}} &= [A_\la]_q.\\
  \end{align*}
  The statements about $B_\la$ are immediate consequences of the fact
  that $B_\la$ is isomorphic to the twist of $A_\la$ by $\la$.
\end{proof}

\begin{remark}\label{rem:one-is-enough}
  This proposition can be used to compute $P_{A_\la,n}$, etc., for all $n$,
  since the number of points of an elliptic curve in an extension of a finite
  field is determined by that of the ground field and $\la \in \F_q$ is a
  square in $\F_{q^n}$ if and only if $\la$ is a square in $\F_q$ or $n$
  is even.
\end{remark}
  
We complete these results by proving analogous ones for the case $\la = -1$.
Neither
the equation $y^2 = x^3 - 2x^2 + \frac{\la}{\la+1}x$ nor any twist gives
an elliptic curve in this case, but we can still describe varieties isogenous
to the Kummer surface by similar formulas.

\begin{definition}\label{def:kmin1}
  Let $K_{-1}, L_{-1}$ be the surfaces defined by
  $$ v^2 = z_0z_1z_2(-z_0+z_1)(z_1+z_2)(z_0+z_2), \ 
  v^2 = -z_0z_1z_2(z_0+z_1)(-z_0+z_2)(z_1+z_2),$$
  and $A_{-1}, B_{-1}$ the affine patches $z_0 \ne 0$.
  Let $a_{-1,p} = p+1-[E]_p$, where $E$ is the elliptic curve with affine
  equation $y^2 = x^3 - x$.
\end{definition}

As before, by exchanging $z_1, z_2$ we see that $K_{-1}, L_{-1}$ are quadratic
twists of each other by $-1$.  On the other hand, the map $(v:-z_0:z_1:z_2)$
is an isomorphism $K_{-1} \to L_{-1}$.

\begin{proposition}\label{prop:count-m1}
  For all odd primes $p$ we have
  $[K_{-1}]_p = [L_{-1}]_p = p^2 - \phi(-1) p + 1 + a_{-1,p}^2$.
\end{proposition}

\begin{proof}
  For $p \equiv 3 \bmod 4$, we have stated above that $K_{-1}, L_{-1}$ are
  isomorphic to their twists by $-1$, which is not a square in $\F_p$, so
  the number of points is $p^2+p+1$.  This is as claimed, since $a_{-1,p} = 0$
  for such $p$.
  
  In the case $p \equiv 1 \bmod 4$, the
  argument is very similar to that given to prove Theorem~\ref{thm:count-kl}.
  The two surfaces are
  isomorphic, so we only consider $K_{-1}$.  Again we begin with the
  fibration $(z_0:z_1)$, for which the general fibre is defined by
  \begin{equation}
    y^2 = x^3 + (-t^3+t)x^2 + (t^5-2t^4+t^3)x
  \end{equation}
  and there are
  three fibres of type $I_2^*$ (since the $I_0^*$ and $I_2$ of the generic
  case come together).  All components of the singular fibres
  are rational.  We consider the quotient of this elliptic surface by
  $(0:0)$, obtaining a surface defined by
  \begin{equation}
    y^2 = x^3 + (2t^3-2t))x^2 + (t(t-1)^2)^2 x.
  \end{equation}
  This surface has an $I_4^*$ fibre at $1$ and $I_1^*$ at $0, \infty$; the action
  of Galois on all three is trivial. % (here we use that $p \equiv 1 \bmod 4$).
  Let $\widetilde S_{-1}$ be the K3 surface given
  by the minimal desingularization of this surface.  Then as in the proof of
  Theorem~\ref{thm:count-kl} we have
  $[\widetilde K_{-1}]_p = [\widetilde S_{-1}]_p$.
  
  On the other hand, we consider the Kummer surface $\Kum_{-1}$ of
  $E_{-1} \times E_{-1}$, where $E_{-1}$ is defined by $y^2 = x^3 - x$.  It can
  be defined in $\P(3,1,1,1)$ by $v^2 = (z_0^3-z_0z_2^2)(z_1^3-z_1z_2^2)$.
  The map defined by
  \begin{equation}
    \begin{split}
&((z_0+z_2)(4z_0^2z_1 + z_0z_1^2 - z_1^2z_2 - z_0z_2^2 - 
      2z_1z_2^2 - z_2^3)/4:\\
      &\quad (z_0z_1 - z_0z_2 - z_1z_2 - z_2^2)(3z_0z_1 - z_0z_2 - z_1z_2 - z_2^2))/3\\
    \end{split}
  \end{equation}
  induces an elliptic fibration on the minimal desingularization whose general
  fibre is isomorphic to that above.  Thus $[\dsKum_{-1}]_p = [\widetilde S_{-1}]_p$.
  
  But as before $[\dsKum_{-1}]_p = p^2 + 18p + 1 + a_{-1,p}^2$.  
  The singular subscheme of $K_{-1}$ has degree $19$ and all the exceptional
  curves are defined over~$\F_p$, so
  $$[K_{-1}]_p = [\widetilde K_{-1}]_p -19p = [\dsKum_{-1}]_p - 19p = p^2 - p + 1 + a_{-1,p}^2$$
  as claimed.
\end{proof}

We now demonstrate a relation between the $a_{\la,p}$ and the coefficients
of an eigenform.  This will be used in Section~\ref{sec:first} to convert an
expression for the number of points derived from a fibration to one in
terms of the coefficients.
\begin{proposition}\label{prop:finish-proof}
$p+\sum_{\la=1}^{p-2} \phi(\la) (a_{\la,p}^2-p) = -b_p$, where
$b_p$ is the eigenvalue of $T_p$ on the newform of weight~$4$ and
level~$8$.
\end{proposition}
 
\begin{proof}
  Consider the total space $T_L$ of the family of
  $L_\la$ (Definition~\ref{def:kla-lla}) as a double cover of
  $\P^2 \times \P^1$.  On the one hand, this is fibred
  by $\la$, with the fibre at $\la$ having $p^2+p+1+\phi(\la)(a_{\la,p}^2-p)$
  points for $\la \ne 0, -1, \infty$ by Corollary~\ref{cor:count-ll};
  it is easily seen that the fibres
  at $0, -1, \infty$ all have $p^2+p+1$ points.  (We do not obtain the
  formula of Proposition~\ref{prop:count-m1} for the fibre at $-1$
  because the factor of $\la+1$ makes it a copy of $\P^2$.)

  On the other hand, using the methods of~\cite{c-sz} we can construct
  a resolution and verify that it is a rigid Calabi-Yau threefold, modular
  of level~$8$, and hence obtain the formula $p^3+2p^2+p+1-b_p$ for
  $[T_L]_p$, valid for all $p>2$.
  (The results of~\cite{c-sz} are stated for double covers
  of $\P^3$ for which the components of the branch locus meet transversely
  in smooth loci; however, the arguments apply without change to double
  covers of any smooth threefold and components of the branch locus meeting in
  unions of smooth components that meet transversely.)  Comparing the two sides
  we obtain the equality
  $$(p+1)(p^2+p+1) + \sum_{\la=1}^{p-2} \phi(\la) (a_{\la,p}^2-p) = p^3+2p^2+p+1-b_p;$$
  now subtract $p^3+2p^2+p+1$ from both sides.
\end{proof}

\begin{remark}\label{rem:comes-from-fop}
  Perhaps this result can be deduced from those of~\cite{fop}, but
  this type of geometric argument can be applied in many situations where
  the methods of~\cite{fop} cannot.
\end{remark}
  
The K3 surfaces considered so far in this section are important for the
calculations in Section~\ref{sec:first}.  In the remainder of this section,
we make some small modifications in order to study K3 surfaces that arise
in Section~\ref{sec:32}.

\begin{definition}\label{def:kla-lla-new} Let $M = M_\la$ and $N_\la$ be the
  surfaces in $\P(3,1,1,1)$  defined by
  $$\begin{aligned}
    M_\la&: t^2 = xyz(x+y)(y+z)(-x+z),\\
    N_\la&: t^2 = \la x(-x+\la z)y(-y+z)(x+2y-z)(-x-2y+(\la+1)z).
  \end{aligned}
  $$
\end{definition}
% this matches k3qss[1] in level32 by means of the automorphism of P^2
% defined by (x-2y+z:x:z)

%% Now, $M = M_{\la}$ is independent of $\la$ and
%% $(M_\la \times N_\la)/\sigma \cong (k_\la \times \ell_\la)/\sigma$.

\begin{proposition}\label{prop:count-k} $[M]_p = p^2+p+1+a_{3,p}$, and the
  branch locus of $M$ has $6p-5$ points over~$\F_p$.
\end{proposition}

\begin{proof}
  Observe that $M$ is the same surface
  as $M_{-1}$ (Definition~\ref{def:kmin1}) up to a change of variables.  We showed
  in Proposition~\ref{prop:count-m1} that $[M]_p = p^2 - \phi(-1) p + 1 + a_{-1,p}^2$.
  Since $y^2 = x^3 - x$ is the unique elliptic curve of conductor $32$
  up to isogeny and $a_{-1,p}$ is the trace of Frobenius at $p$ for this
  curve, we see that $a_{-1,p} = a_{2,p}$.  In light of Lemma~\ref{lem:coef-mf},
  this implies our claim for $p \equiv 1 \bmod 4$.  For $p \equiv 2, 3 \bmod 4$,
  both sides are equal to $p^2+p+1$.
\end{proof}

Recall the notation (Definition~\ref{def:count-double-cover})
$P_{V,i}, Z_{V,i}, N_{V,i}$ for the number
of points on the base of a double cover over $\F_{p^i}$ that pull back
to $2, 1, 0$ points on the cover respectively.
In these terms,
we may rephrase Proposition~\ref{prop:count-k} as
saying that $Z_{M,1} = 6p-5, P_{M,1} = (p^2-5p-4+a_{3,p})/2, N_{M,1} = (p^2-5p-4-a_{3,p})/2$.
We also note that $Z_{N_\la,1} = p^2+p+1$ for $\la = 0$, while it is
$6p-7$ for $\la = \pm 1$ and $6p-9$ for other values of $\lambda$ (the
difference is that two sets of three lines in the branch locus are concurrent
for $\pm 1$ but not on other fibres).

The Picard number
of $N_\la$ is generically $19$: this can be seen by constructing an elliptic
fibration on it with two bad fibres each of type $I_0^*, I_4, I_2$, and a
section of infinite order.
In fact, the Picard lattice of $N_\la$ is a sublattice of index
$2$ of that of the Kummer surface of the square of an elliptic curve
without complex multiplication.
%\begin{remark}\label{rem:lla-is-twist}
The surface $N_\la$ is isogenous to the Kummer surface of
$E_\la \times E_\la^\tau$, where $E_\la$ is an elliptic curve with
$j$-invariant
$(-4\la^2+16)^3/\la^4$ and $\tau$ is its quadratic twist by $-\la^3+\la$. 
  %% This observation leads to an unnecessarily
  %% complicated method of counting the $\F_p$-points on $V_{32}$,
  %% but it also suggests an interesting equality that will be discussed
  %% in a future paper.
%\end{remark}

It is not so easy to give a useful formula for $[N_\la]_p$.
To avoid having to do so,
we consider the total space $\L$ of the family of the $N_\la$ inside
$\P(3,1,1,1) \times \P^1$.  In other words, we regard $\la$ as the ratio of
the two coordinates of $\P^1$.
Let the coordinates on this space be $t,z_0,z_1,z_2,u,v$: then $\L$ is
defined by the equation
$$t^2v^3 = uz_0(-vz_0+uz_2)z_1(-z_1+z_2)(z_0+2z_1-z_2)(-vz_0-2vz_1+(u+v)z_2).$$
We define a rational map $\rho: \L \dashrightarrow \P^1$ by $(2z_1-z_2:z_2)$.
It is easily checked that the base scheme consists of two rational
curves that meet in a single point; it thus has $2p+1$ points mod
$p$ for all $p$.  As in Proposition~\ref{prop:finish-proof},
we will use a different fibration on $\L$ to count its $\F_p$-points.

\begin{proposition}\label{prop:fibres-rho}
  Let $x \in \F_p$ with $x^3 - x \ne 0$
  and let $\rho_x$ be the fibre of $\rho$ at $(x:1)$.  Then
  $\rho_x$ has $p^2+4p+1+\phi(x^3-x) a_{3,p}$ points.
  The fibres of $\rho$ at $0, \pm 1, \infty$ have
  $p^2 + 3p + 1, 2p^2+2p+1, 2p^2+2p+1$ points respectively.
\end{proposition}

\begin{proof}
  We consider the affine patch of $\rho_x$ where $z_2, v$ are nonzero.
  We may view this patch of $\rho_x$
  as being inside $\A^3$, which in turn we think of as the affine patch
  of $\P(3,1,1,1)$ where the last coordinate is nonzero.  The projective
  closure is defined by the equation
  $$t^2 = \frac{-x^2+1}{4}z_0z_1z_2(z_0+xz_2)(z_0-z_1)(z_0-z_1+xz_2).$$
  Replacing $z_2, t$ by $z_2/x, t/2x$ and multiplying through by $(x/2)^2$
  converts this to
  \begin{equation}\label{threefold-e32}
    t^2 = (-x^3+x)z_0z_1z_2(z_0+z_2)(z_0-z_1)(z_0-z_1+z_2),
  \end{equation}
  and replacing $z_0-z_1$ by $z_1$ converts $z_1$ to $z_0-z_1$ and hence
  changes this to the equation for $K$, twisted by $-x^3+x$,
  up to the order of variables.
  
  Note further that if $p \equiv 1 \bmod 4$ then
  $\phi(x^3-x) = \phi(-x^3+x)$, and if $p \equiv 3 \bmod 4$ then
  $a_{3,p} = 0$, so we may replace $\phi(-x^3+x)$ by $\phi(x^3-x)$ in
  both cases.  With these observations, the proof for the general fibres
  reduces to routine bookkeeping.
  
  As for the bad fibres of $\rho$, the fibre at $1$ consists of two
  components, one supported at $t = z_1 - z_2 = 0$ and one at
  $z_1 - z_2 = v = 0$.  The total number of points is
  $p^2+2p+1+p^2+p+1-(p+1) = 2p^2 + 2p + 1$; similarly for $\rho_{-1}$
  with $z_1 - z_2$ changed to $z_1$.
  
  The fibre at $\infty$ has components at
  $t^2v + (z_0z_1(z_0+z_1))^2u = z_2 = 0$ and $z_2 = v = 0$.  To count the
  points on the first of these, note that for fixed $t, z_0, z_1$ and
  $z_2 = 0$ we get one solution for $u, v$ if $t \ne 0$ or
  $z_0z_1(z_0+z_1) = 0$ and $p+1$ otherwise.  Thus the total is
  $p^2+p+1+3p = p^2 + 4p + 1$.  The two components intersect along
  $z_2 = v = z_0z_1(z_0+z_1)$, so they share $3p+1$ points.  The
  second component has $p^2+p+1$ points, so the total is $2p^2 + 2p + 1$.
  
  Finally, the fibre at $0$ is defined in $\P(3,1,1,1) \times \P^1$ by
  $$4t^2v^3 = z_0^2 z_2^2(z_0^2uv^2 + 2z_0z_2uv^2 + z_2^2(u^3+2u^2v)), \quad 2z_1 - z_2 = 0.$$
  We count its points with the help of the projection to $\P(3,1,1,1)$.
  It is readily checked that the fibre at $0$ is supported on a smooth
  rational
  curve, giving $p+1$ points, and that the fibre at $\infty$ is supported on
  two smooth rational curves that meet at $(1:0:0:0)$, giving $2p+1$ points.
  The fibre at $(\alpha:1)$ consists of two rational curves that meet at
  $3$ rational points.  The components are defined over $\F_p(\sqrt{\alpha})$,
  so we find $p+1+(p-2)\phi(\alpha)$ points.  The map has empty base scheme,
  and $\phi(\alpha)$ is $+1$ and $-1$ equally often, so the total number of
  points is $p(p+1)+2p+1 = p^2 + 3p + 1$ as claimed.
\end{proof}

\begin{corollary}\label{cor:totl-is-mtimese-mod-sigma}
  Let $E$ be the elliptic curve $y^2 = x^3 - x$ with the involution
  $(x,y) \to (x,-y)$.  Then $\L$ and $(M \times E)/\sigma$ are birationally
  equivalent.
\end{corollary}

\begin{proof}
    This follows from the
  equation (\ref{threefold-e32}) for the fibre of $\rho_x$.
\end{proof}

\begin{corollary}\label{cor:count-script-l}
  The total space $\L$ has $p^3+6p^2-3p+1-a_{4,p}-pa_{2,p}$ points over~$\F_p$.
\end{corollary}

\begin{proof} By the proposition, there are $(p-3-a_{2,p})/2$ fibres with
  $p^2+4p+1+a_{3,p}$ points and $(p-3+a_{2,p})/2$ with $p^2+4p+1-a_{3,p}$ points,
  in addition to the points of the bad fibres.  The base scheme having
  $2p+1$, the total number of points is then
  $$\begin{aligned}
    &\frac{1}{2} ((p-3-a_{2,p})(p^2+4p+1+a_{3,p}) + (p-3+a_{2,p})(p^2+4p+1-a_{3,p}))\cr
    &\quad + p^2+3p+1+3(2p^2+2p+1)-p(2p+1).\cr
  \end{aligned}$$
  Simplifying and applying Lemma~\ref{lem:coef-mf} gives this result.
\end{proof}

\section{The first example: level $8$}\label{sec:first}
  Our first example of a modular fivefold of level~$8$ is the double cover
$V_8$ of $\P^5$ defined by the equation
$$t^2 = \prod_{i=0}^5 x_i (x_i + x_{i+1}).$$
This equation is reminiscent of the first arrangement of $8$
hyperplanes in $\P^3$ given in the table of~\cite[{p.{} 68}]{meyer}, which
likewise defines a variety that is modular of level~$8$.  However, the
pattern does not persist beyond dimension $5$.

\begin{remark}\label{rem:aut-f1}
  The group of automorphisms of the double cover $V_8 \to \P^5$
  is of order $24$: it is generated by the maps $V_8 \to V_8$
  taking $(t:x_0:\dots:x_5)$ to
    $$(-t:x_0:\dots:x_5), \quad (t:x_1:\dots:x_5:x_0), \quad
    (t:x_5:x_4:\dots:x_0).$$
  This is checked by verifying that the only automorphisms of the dual
  $\check \P^5$ preserving the $12$ points corresponding to the $12$
  hyperplanes are the obvious ones.
\end{remark}

We will prove its modularity in the following precise form.
\begin{definition}\label{def:newforms}
  Let $a_p, b_p$ be the Hecke eigenvalues for the unique newforms of level
  $8$ and weight~$6, 4$ respectively.  (In the LMFDB \cite{lmfdb} these
  are {\tt 8.6.a.a} and {\tt 8.4.a.a}.)
\end{definition}

\begin{theorem}\label{thm:main-first}
  \begin{equation}\label{f1-formula}
    [V_8]_p = p^5 + p^4 + p^3 + p^2 + p + 1 - a_p - (b_p + \phi(-1) p) p.
  \end{equation}
\end{theorem}

\begin{remark} This is related to a recent result of
  Li, Long, and Tu \cite[Theorem 5, first row of table]{llt}.
  They identify a Galois
  representation for which the trace of Frobenius is given by the same
  formula $a_p + (b_p + \phi(-1)) p$.  The methods of proof are also
  related, since their expression of a value of $_6F_5$ in terms of a
  sum over products of $_3F_2$ can be interpreted in terms of the
  fibration that we will discuss in the next section.
\end{remark}

\subsection{Proof of modularity}\label{subsec:proof-mod-f1}
We will prove Theorem~\ref{thm:main-first} by means of a fibration by
quotients of products of two K3
surfaces.  By means of results of Section~\ref{sec:geom-points},
we will relate these to the Kummer surface of the square of an
elliptic curve and use this to express the number of points in terms of
hypergeometric functions over a finite field, thus reducing to results
of~\cite{fop}, especially Theorem 1.1.
The statements in this section have been proved; but, since it is easy to
make mistakes in such things, we also verify most of the statements
of this section numerically for small $p$ in the file
%\ref{prop:count-prod},~\ref{prop:count-kum},~\ref{prop:fib-on-kum},
%\ref{prop:fib-on-k},~\ref{thm:count-kl},~\ref{thm:count-ll},
%\ref{cor:count-f},~\ref{prop:count-fm1
{\tt code-level-8.mag} in \cite{magma-scripts} by counting points directly
and with the aid of the various fibrations and other geometric constructions.

\begin{definition} Let $\pi$ be the rational map $V_8 \dashrightarrow \P^1$
  defined by $(x_0:x_3)$.
\end{definition}

The justification for this definition is that if we set
$x_0 = \lambda, x_3 = 1$
in the equations of the $12$ hyperplanes along which the double cover
$V_8 \dashrightarrow \P^5$ is branched, we can divide them into
two sets of $6$ linear forms, one depending only on $\lambda, x_1, x_2$,
the other only on $\lambda, x_4, x_5$.
  
\begin{definition}\label{def:fla}
  Let $F_\la$ be the affine patch $x_3 \ne 0$ of the fibre of $\pi$ above
  $(\la:1)$.
\end{definition}

\begin{proposition}\label{prop:fib-is-prod} For $\la \in \F_p$ with $\la \ne 0, -1$,
  the fibre of $\pi$ at $(\la:1)$ is birationally equivalent
  to $(K_\la \times L_\la)/\sigma$, where $\sigma$ is the automorphism 
  $$((v:z_0:z_1:z_2),(w:y_0:y_1:y_2)) \to ((-v:z_0:z_1:z_2),(-w:y_0:y_1:y_2)).$$
\end{proposition}

\begin{proof} As above, if we substitute $x_0 = \la, x_3 = 1$ in the equations
  for $V_8$, the variables separate and the $12$ linear forms can be expressed
  as $6$ in $x_1, x_2$ and $6$ in $x_4, x_5$.  Thus there is an obvious map
  $\tau$ of degree $2$ from the product of the two double covers of $\A^2$
  branched along these loci to the fibre of $\pi$, and
  $\tau \circ \sigma = \tau$.  The result follows by taking the projective
  closure.
\end{proof}

\begin{proposition}\label{prop:count-prod}
  For all odd primes $p$ and all $\la \ne 0, -1 \in \F_p$ we have
  $[F_\la]_p = [A_\la]_p[B_\la]_p - p^2[A_\la]_p - p^2[B_\la]_p + 2p^4$.
\end{proposition}

\begin{proof}
  Apply Lemma~\ref{lem:count-q-dc} to  Proposition~\ref{prop:fib-is-prod}.
\end{proof}

  %%% why?
  %% As an alternative we show how to prove it by means of character sums.  We have
  %% \begin{equation}
  %% \label{char-sums}
  %% \begin{split}
  %%   [F_\la]_p &= \sum_{y_1, y_2, y_4, y_5 \in \F_p} 1+\phi(\la(\la+1) y_1y_2(\la+y_1)(y_1+y_2)(y_2+1) \times\\
  %%   &\qquad\qquad\qquad\qquad (\la+1)y_4y_5(1+y_4)(y_4+y_5)(y_5+\la)),\\
  %%  [A_\la]_p &= \sum_{y_1,y_2 \in \F_p} 1 + \phi((\la+1) y_1y_2(\la+y_1)(y_1+y_2)(y_2+1)),\\
  %%  [B_\la]_p &= \sum_{y_4,y_5 \in \F_p} 1 + \phi(\la(\la+1)y_4y_5(1+y_4)(y_4+y_5)(y_5+\la)).\\
  %% \end{split}
  %% \end{equation}
  %%  Thus
  %%  \begin{equation}
  %%  \label{f-ab-decomp}
  %%  \begin{split}
  %%  [F_\la]_p &= [A_\la]_p[B_\la]_p 
  %% - \sum_{y_1,y_2 \in \F_p}\phi((\la+1) y_1y_2(\la+y_1)(y_1+y_2)(y_2+1))\\
  %% &\qquad - \sum_{y_4,y_5 \in \F_p}\phi(\la(\la+1)y_4y_5(1+y_4)(y_4+y_5)(y_5+\la))\\
  %% &= [A_\la]_p[B_\la]_p - p^2([A_\la]_p-p^2) - p^2([B_\la]_p-p^2)\\ \end{split}
  %%  \end{equation}
  %% as claimed.

Combining Corollary~\ref{cor:count-ll} with Proposition~\ref{prop:count-prod}
and Lemma~\ref{lem:ka-lb} (together with the 
observation that $[F_0]_p = p^4$, because $x_0=0$ is one of the hyperplanes and
so there is one $\F_p$-point of $F_0$
for each point of $\A^4(\F_p)$), we immediately obtain
the following:

\begin{corollary}\label{cor:count-f}
  Let $\la \ne -1$.  Then $F_\la$
  has $p^4 + \phi(\la) (a_{\la,p}^2-p)^2$ points.
\end{corollary}

In the case $\la = -1$,
the separation of the variables still allows us to write
$[F_\la]_p$ in terms of a K3 surface and its quadratic twist by $\la$.
The only difference is that we need to use Proposition~\ref{prop:count-m1}
instead of Lemma~\ref{lem:ka-lb} and Corollary~\ref{cor:count-ll}.

\begin{proposition}\label{prop:count-fm1}
  On $F_{-1}$ there are
  $p^4 + (2p-a_{-1,p}^2)^2$ points for $p \equiv 1 \bmod 4$
  and $p^4$ for $p \equiv 3 \bmod 4$.
\end{proposition}

\begin{proof}
  The count of points on $F_{-1}$ follows by applying
  Proposition~\ref{prop:count-prod}
  to the result of Proposition~\ref{prop:count-m1}
  with the help of Lemma~\ref{lem:ka-lb} (which clearly applies in the case
  $\la = -1$ now that the necessary surfaces have been defined),
  as in Corollary~\ref{cor:count-f}.
  %% Alternatively, we may use
  %%   Proposition~\ref{prop:fib-is-prod} together with Lemma~\ref{lem:count-q-dc}.
  %%   Let $k_+, k_0, k_-$ be the number
  %%   of points of the affine patch $z_0 = 0$ of $\P^2$ where the branch function
  %%   of $K_{-1}$ is a nonzero square, $0$, or a nonsquare respectively: then
  %%   $[F_{-1}]_p = p^4 + (k_+-k_-)^2.$  The first
  %%   statement of the proposition says that $k_+ - k_- = 2p - a_{-1,p}^2$ and the
  %%   result follows.
\end{proof} 

We are now ready to prove Theorem~\ref{thm:main-first}.  For simplicity we will
only write out the proof in the case of $p \equiv 3 \bmod 4$; the case
$p \equiv 1 \bmod 4$ is very similar but requires slightly more work
to keep track of the $\la = -1$ terms.

\begin{proof}
  With $p \equiv 3 \bmod 4$, the formula
in Theorem~\ref{thm:main-first} becomes
$$\sum_{i=0}^5 p^i - a_p - pb_p + p^2$$ (recall that $a_p, b_p$ are the Hecke
eigenvalues for the newforms of weight~$8$ and level~$6, 4$ respectively).
By \cite[Theorem 1.1]{fop}, for $p \equiv 3 \bmod 4$ we have
$$a_p = -p^4 \sum_{\la=2}^{p-1}\phi(-\la)\fthreetwo(\la)^2 + p^2 - pb_p,$$ so we need
to show that
\begin{equation}\label{eqn:f1-formula-fthreetwo}
  [V_8]_p = \sum_{i=0}^5 p^i - p^4 \sum_{\la=2}^{p-1}\phi(-\la)\fthreetwo(\la)^2.
\end{equation}

We count the points on $V_8$ by means of the fibration $\pi$.
The hyperplane $x_3 = 0$ has
$p^4+p^3+p^2+p+1$ points, in bijection with those of the hyperplane $x_3 = 0$
in $\P^4$.  The affine patch $x_3 \ne 0$ of the fibre at $0$ has $p^4$ points.
So, by Corollary~\ref{cor:count-f} and Proposition~\ref{prop:count-fm1},
the total number of points is
$$\sum_{i=0}^4 p^i + 2p^4 + \sum_{\la = 1}^{p-2} p^4 + \phi(\la)(a_{\la,p}^2-p)^2.$$
  Combining the $p^4$ terms, changing $\la$ to $-\la$ in the summation, and
  using (\ref{f-a}) shows that this is the same as the right-hand side of
  (\ref{eqn:f1-formula-fthreetwo}).
  This completes the proof of Theorem~\ref{thm:main-first} in the case
  $p \equiv 3 \bmod 4$.
  %% As already mentioned, the proof for
  %% $p \equiv 1 \bmod 4$ is very similar, requiring the use of 
  %% the familiar fact \cite[Theorem 18.5]{ireland-rosen}
  %% that $a_{-1,p} = \pm 2a$, where as
  %% before $a$ is the positive odd integer such that $p - a^2$ is a square.
\end{proof}

\subsection{Constructing an apparent rigid Calabi-Yau fivefold from $V_8$}\label{subsec:rigid-f1}
Following a method of Burek \cite{burek}, we will construct a
fivefold realizing the same newform of weight~$6$
and level~$8$ as a quotient of $V_8$.  Recall that there is an action of
$D_6$ on $\P^5$ that preserves the set of components of the branch locus of the
map $V_8 \to \P^5$.
We will consider the
quotients $Q_1, Q_2, Q_3$ of $V_8$ by representatives
$\iota_1, \iota_2, \iota_3$
of each of the three conjugacy classes of involutions in $D_6$.
All of these involutions commute with the map that
exchanges the sheets of the double cover, so the quotients are still double
covers of a quotient of $\P^5$: let these quotients be $R_1, R_2, R_3$.

By computing $[Q_1]_p, [Q_2]_p$ for small $p$ we will obtain information
about the action of the $\iota$ on the cohomology of a resolution, leading
us to believe that $Q_3$ realizes the newform in question.
In Theorem~\ref{thm:q3} we will prove a formula expressing $[Q_3]_p$
in terms of powers of $p$, the class of $p \bmod 4$, and the Fourier
coefficient $a_p$ of the newform.
In the file {\tt quotient-level8.mag} in \cite{magma-scripts} we
verify the statements 
of Proposition~\ref{prop:count-q-r-small-p}, Conjecture~\ref{conj:q2},
Theorem~\ref{thm:q3} for small $p$ by comparing the formulas we give
to the number of points obtained directly by counting.

\begin{remark}\label{rem:prob-cy}
  We expect that $Q_3$ has a resolution which is a strongly rigid
  Calabi-Yau fivefold, but are unable to prove this.
  %% This is discussed further in Remark~\ref{rem:q3-hard}.
\end{remark}

The differential $D$ on $\P^5$ given by
\begin{equation}\label{hol-diff-p5}
  (-1)^j\frac{x_j^5}{\prod_{\substack{i=0\\i \ne j}}^5 x_i} \bigwedge_{\substack{i=0\\i \ne j}}^5 d(x_i/x_j)
\end{equation}
(\cite[Remark III.7.1.1]{h}) is independent of $j$ and has divisor
$-\sum H_i$, where $H_i$ is the hyperplane $x_i = 0$.
Pulling it back to $V_8$, we get a differential whose divisor is 
$R - 2\sum H_i$, where $R$ is the ramification locus.  This is the divisor
of $t/\prod_{i=0}^5 x_i$, so we obtain a differential
\begin{equation}\label{hol-diff-f1}
  (-1)^j\frac{x_j^6}{t} \bigwedge_{\substack{i=0\\i \ne j}}^5 d(x_i/x_j)
\end{equation}
on $V_8$ whose divisor is $0$.
If $V_8$ admits a crepant resolution $\tilde V_8$, then $D$ pulls back to
a differential on $\tilde V_8$ with divisor $0$, and so we can compute
the action of $\iota_i$ on $H^{5,0}(\tilde V_8)$ from its action on $D$.
Since we do not know this (indeed, we believe that it is false, for reasons
explained in Section~\ref{subsec:geom-cy-res}), this
computation can be used only as motivation, not in any proofs.

To see that $\iota_1^*(D) = -D$, note that
$t/\prod_{i=0}^5 x_i$ is invariant under $\iota_1$, while
exchanging two variables changes the sign of $D$ (choose
an expression for $D$ in which $x_j$ is not one of the two variables).
Since $\iota_1$ gives
an odd permutation, it acts as $-1$ on the pullback of $D$ to $V_8$ and
we do not expect to see
the form of weight~$6$ in the cohomology of the quotient.

First, we examine the central element $\iota_1: x_i \to x_{i+3}$.
The invariant ring for the action of $\iota_1$ on $\P^5$ is generated by the
polynomials $x_i+x_{i+3}, x_i^2 + x_{i+3}^2$ for $0 \le i \le 2$
and $x_ix_j + x_{i+3}x_{j+3}$ for
$0 \le i < j \le 2$.  Thus we may view the quotient map $\P^5 \to R_1$
as a map to a subvariety of $\P(1,1,1,2,2,2,2,2,2)$.  None of the $12$
hyperplanes is fixed by the involution, so we obtain $6$ branch divisors,
each of which is defined by an equation of degree~$2$.

\begin{proposition}\label{prop:count-q-r-small-p}
  For $p$ an odd prime less than $20$, both
  $Q_1$ and $R_1$ have $\sum_{i=0}^5 p^i$ points over~$\F_p$.
\end{proposition}

\begin{proof} First find single polynomials $s_i$ defining each
  of the branch components as a subvariety of $R_1$; then, for each $p$,
  enumerate the $\F_p$-points of $R_1$, evaluate the product of the $s_i$ at
  each, and sum the Kronecker symbols to obtain the point count of $Q_1$.
  All of this is easily done in Magma when $p$ is small.
\end{proof}

\begin{remark}\label{rem:proof-boring}
  It would be tedious but not difficult to prove this proposition
  for all odd primes $p$.  There being no lack of varieties already known
  to have $\sum_{i=0}^5 p^i$ points over~$\F_p$ for all $p$,
  such a result would be of little interest, and so we
  have not done this.
\end{remark}
%% In Remark~\ref{resol-cohom} we suggested that $V_8$ has a
%% resolution with $h^{5,0} = h^{4,1} = h^{1,4} = h^{0,5} = 1$; if so, presumably
%% $\iota_1$ acts as $-1$ on $H^5$ of this resolution.

Next we consider $\iota_2: x_i \to x_{3-i}$.  As an involution of $\P^5$
this is conjugate to $\iota_1$, and again it gives an odd permutation of the
variables and hence acts as $-1$ on $D$,
but it is different as an automorphism of
$V_8$.  Indeed, it fixes $2$ of the $12$ components of the branch locus,
so we get $7$ branch divisors, of which $5$ have degree $2$ and $2$ have
degree $1$.  Again we find that $[R_2]_p = \sum_{i=0}^5 p^i$, but this time
$[Q_2]_p = \sum_{i=0}^5 p^i - pb_p$ for $p < 20$, where as before $b_p$ is
the Hecke eigenvalue for the newform of weight~$4$ and level~$8$.  This
suggests that $\iota_2$ acts as $+1$ on $H^{4,1} \oplus H^{1,4}$ and as
$-1$ on $H^{5,0} \oplus H^{0,5}$.
\begin{conjecture}\label{conj:q2}
  $[Q_2]_p = \sum_{i=0}^5 p^i - pb_p$ for all odd $p$.
\end{conjecture}

\begin{remark}\label{rem:why-not-q1-q2} Again, the methods that prove
  Theorem~\ref{thm:q3} could be applied to prove Conjecture~\ref{conj:q2},
  but this seems much less novel than Theorem~\ref{thm:q3}
  and so we have not carried it out.
\end{remark}

Let $\iota_3 = \iota_1 \iota_2$.  In view of
Proposition~\ref{prop:count-q-r-small-p}
and Conjecture~\ref{conj:q2}, we expect that $\iota_3$
acts as $+1$ on $H^{5,0} \oplus H^{0,5}$ of a resolution of $V_8$
and as $-1$ on $H^{4,1} \oplus H^{1,4}$.
This time the ring of invariants is generated by
$x_0, x_1+x_5, x_2+x_4, x_3, x_1^2+x_5^2, x_2^2+x_4^2,x_1x_2+x_4x_5$, so the
quotient map from $\P^5$ goes to $\P(1,1,1,1,2,2,2)$; the image is 
a hypersurface $H$ of degree $4$.  Two of the $12$ components of
the branch locus are fixed by $\iota_3$ and map to divisors in $H$ cut out
by equations of degree $1$; the other $10$ are exchanged in pairs and map
to divisors defined by equations of degree $2$.  Numerically this suggests
a singular Calabi-Yau variety: the canonical divisor of $H$ would be
$\O(-4 \cdot 1 - 3 \cdot 2 + 4) = \O(-6)$, while the branch locus has class
$\O(12) = -2K_H$, but in view of the large
singular locus of $H$ and the branch divisors
this is not a proof.  By directly counting points
we obtain the following formula, which we will prove in the next section.
\begin{theorem}\label{thm:q3}
  %$[R_3]_p = \sum_{i=0}^5 p^i$ and
  $[Q_3]_p = \sum_{i=0}^5 p^i - a_p - \phi(-1)p^2$ for all odd $p$,
  where as before $a_p$ is the eigenvalue of $T_p$ on the newform of
  weight~$6$ and level~$8$.  
\end{theorem}

This and 
Propositions~\ref{prop:count-q-la-res},
\ref{prop:count-q-la-nonres}, \ref{prop:q3-minus-f1} are verified numerically in
the file {\tt count-quotient.mag} in \cite{magma-scripts} by comparing the
formulas we give to point counts computed directly from equations for the
varieties.
%% \begin{remark}\label{rem:why-not-geom}
%%   In principle this is a theorem about the Galois action on the \'etale
%%   cohomology of a resolution of $Q_3$.  We do not prove it in this way,
%%   because that would require a very detailed understanding of such a
%%   resolution, which can only be constructed after many blowings-up.
%% \end{remark}
\subsection{Proof of Theorem~\ref{thm:q3}}\label{subsec:proof-thm-q3}

Since $\iota_3$ preserves the fibration $(x_0:x_3)$,
the natural approach is to study the fibres of the
map $\pi_3: Q_3 \to \P^1$ given by $(x_0:x_3)$.  These are in fact the
quotients of the fibres of $\pi$ by $\iota_3$.  Let $Q_{3,\la}$ be the
fibre of $Q_3$ above $(\la:1)$.  We will count the points of $Q_{3,\la}$ in
terms of $P_{A_\la,1}$, etc.  In particular we show:
\begin{proposition}\label{prop:count-q-la-res}
  Let $\la \ne 0, -1$, and suppose that $\phi(p) = 1$.  Then
  \begin{align*}
    [Q_{3,\la}]_p &= P_{A_\la,1}^2 + P_{A_\la,1} + N_{A_\la,1}^2 + N_{A_\la,1} +    (P_{A_\la,2} - P_{A_\la,1} - N_{A_\la,1}) + \\
    &\qquad P_{A_\la,1} Z_{A_\la,1} +    N_{A_\la,1} Z_{A_\la,1} + Z_{A_\la,1}(Z_{A_\la,1}+1)/2 + 5p(p-1)/2.\\
  \end{align*}
\end{proposition}

\begin{remark} The $P_{A_\la,1}$, etc., are determined in
  Proposition~\ref{prop:pzn-a-b}, while Remark~\ref{rem:one-is-enough}
  explains how to use that proposition to evaluate $P_{A,\la,2}$.
\end{remark}

\begin{proof} This comes down to considering the pairs of $\F_{p^2}$-points
  on $(A_\la \times B_\la)/\pm 1$ whose images on $Q_{3,\la}$ are rational.
  Let $\mu_\la$ be the $\F_p$-isomorphism $K_\la \to L_\la$ taking
  $(v:z_0:z_1:z_2)$ to $(rw:y_0:y_2:y_1)$, where $r$ is a fixed square root
  of $\lambda$.
  Given a pair of points of $\A^2(\F_p)$ or $\A^2(\F_{p^2})$, we may
  pull them back to $A_\la \times B_\la$ and determine whether their
  inverse images map to rational points of $Q_{3,\la}$.  This occurs in
  the following cases: 
  \begin{enumerate}
  \item A $P$-point $p_A$ of $A$ and a distinct $P$-point $p_B \ne \mu_\la(p_A)$
    of $B$.
    Such a pair gives $2$ points of $Q_{3,\la}$; however, the point
    $(\mu_\la^{-1}(p_B),\mu_\la(p_A)$ of $A_\la \times B_\la$ gives the same points
    of $Q_{3,\la}$,
    so the contribution from this source is $P_{A_\la,1}(P_{A_\la,1}-1)$.
  \item An $N$-point $n_A$ of $A$ and a distinct $N$-point
    $n_B \ne \mu)\la(n_A)$ of $B$.
    Similarly, these give $N_{A_\la,1}(N_{A_\la,1}-1)$ points.
  \item A $P$-point $p_A$ of $A$ and the matching $P$-point $\mu_\la(p_A)$ of $B$.
    Each of these gives $2$ points of $Q_{3,\la}$, for a total of $2P_{A_\la,1}$.
  \item An $N$-point $n_A$ of $A$ and the matching $N$-point $\mu_\la(n_A)$ of
    $B$.  Similarly these give us $2N_{A_\la,1}$.
  \item A $P$-point $p_{A,2}$ of $A$ over $\F_{p^2}$ and $\mu_\la$ of
    its conjugate on $B$.  Since we must exclude the $\F_p$-rational
    points, the number of $P$-points available is
    $P_{A_\la,2} - P_{A_\la,1} - N_{A_\la,1}$ (recall that every element of
    $\F_p$ is a square in $\F_{p^2}$).  Each gives $2$ points of
    $Q_{3,\la}$, but a point and its conjugate give the same points of
    $Q_{3,\la}$, so these give $P_{A_\la,2} - P_{A_\la,1} - N_{A_\la,1}$ points
    of $Q_{3,\la}$.
  \item A $Z$-point of $A$ and a $P$- or $N$-point of $B$.  Because
    a $Z$-point is involved, we are in the branch locus of $Q_{3,\la}$,
    and each such point gives $1$ point of $Q_{3,\la}$.  Thus the number is
    $Z_{A_\la,1}(P_{A_\la,1} + N_{A_\la,1})$.  Note that the images of
    a $P$- or $N$-point of $A$ and a $Z$-point of $B$ are subsumed here,
    the involution $\iota_3$ having the effect of interchanging
    $A_\la$ and $B_\la$.
  \item A $Z$-point $z_A$ of $A$ and a different $Z$-point $z_B \ne \mu_\la(z_A)$
    of $B$.  Again,
    switching these does not change the image, so we obtain a contribution
    of $\binom{Z_{A_\la,1}}{2}$.
  \item A $Z$-point $z_A$ of $A$ and the corresponding $Z$-point $\mu_\la(z_A)$
    of $B$.  There are
    $Z_{A_\la,1}$ possibilities.
  \item A $Z$-point of $A$ whose coordinates generate
    $\F_{p^2}$ and the conjugate of
    its $\mu_\la$-image on $B$.  Every line over~$\F_p$ has $p^2-p$ such
    points and there are $5$ in the branch locus; we must divide by $2$
    since a point and its conjugate map to the same point of $Q_{3,\la}$.
    So the total is $(5p^2-5p)/2$.
  \end{enumerate}
  Adding up the counts from each type gives the result claimed.
\end{proof}

A very similar result holds for $\la$ not a square.
\begin{proposition}\label{prop:count-q-la-nonres}
  Let $\la \ne 0, -1$, and suppose that $\phi(p) = -1$.  Then
  \begin{align*}
    [Q_{3,\la}]_p &= P_{A_\la,1}P_{B_\la,1} + N_{A_\la,1}N_{B_\la,1} +  N_{A_\la,2} + P_{A_\la,1} Z_{A_\la,1} +\\
    &\qquad N_{A_\la,1} Z_{A_\la,1} + Z_{A_\la,1}(Z_{A_\la,1}+1)/2 + 5p(p-1)/2.
  \end{align*}
\end{proposition}

\begin{proof}
  This is identical to the proof of Proposition~\ref{prop:count-q-la-res},
  except for two changes resulting from the fact that $\mu_\la$ is defined
  only over $\F_{p^2}$.  This leads to $P$-points of $B_\la$ over~$\F_p$ matching
  with $N$-points of $A_\la$, rather than $P$-points, and to $N$-points
  of $A_\la$ over $\F_{p^2}$ contributing to the
  count instead of $P$-points.
\end{proof}

Finally we state the results for $\lambda = -1, 0, \infty$.  For
$\la = 0, \infty$, the fibre $Q_{3,\la}$ is contained in the branch
locus and hence has $p^4+p^3+p^2+p+1$ points.
For $\la = -1$, the method is the same as in Propositions~\ref{prop:count-q-la-res},~\ref{prop:count-q-la-nonres}.

\begin{proposition}\label{prop:count-q-la-minus-1}
  We have $[Q_{-1}]_p = p^4 + \phi(\lambda)(a_{\la,p}^2-p)^2 - p^2$.
\end{proposition}

We now have sufficient information to count the $\F_p$-rational points of
$Q_3$.  To avoid duplicating the work already done to prove
Theorem~\ref{thm:main-first}, we start by comparing the point counts on the
fibres of $V_8$ to those on the fibres of $Q_3$.

\begin{proposition}\label{prop:q3-minus-f1}
  For $\la = 0, \infty$ we have $[F_\la]_p = [Q_{3,\la}]_p$, while
  $[F_{-1}]_p - p^2 = [Q_{-1}]_p$.  For other
  values of $\la$, the equality
  $[F_\la]_p - \phi(\la) p(a_{\la,p}^2-p) = [Q_{3,\la}]_p$ holds.
\end{proposition}

\begin{proof}
  Using the results of Corollary~\ref{cor:count-f} and Propositions
 ~\ref{prop:count-q-la-res},~\ref{prop:count-q-la-nonres},
 ~\ref{prop:count-q-la-minus-1}, this is routine.  For example,
  let us suppose that $\la \ne -1$ is a quadratic nonresidue mod $p$.
  In this case we have $[F_\la]_p = p^4 - (a_{\la,p}^2-p)^2$.
  In addition, we have $P_{A_\la,1} = N_{B_\la,1} = (p^2-6p+7+a_{\la,p}^2)/2$, while
  $Z_{A_\la,1} = Z_{B_\la,1} = 5p-7$ and
  $N_{A_\la,1} = P_{B_\la,1} = (p^2-4p+7-a_{\la,p}^2)/2$.  Since $\la$ is a square in
  $\F_{p^2}$, we have $N_{A_\la,2} = N_{B_\la,2} = (p^4 -4p^2+7-a_{\la,p^2}^2)/2$, and
  $a_{\la,p^2} = a_{\la,p}^2 - 2p$ by standard results on elliptic curves
  over finite fields.  Substituting these expressions into
  Proposition~\ref{prop:count-q-la-nonres} shows that
  $$[Q_{3,\la}]_p = p^4 - 2p^2 + 3pa_{\la,p}^2 - a_{\la,p}^4 = [F_\la]_p - p(p-a_{\la,p}^2)$$
  as desired.
\end{proof}

Summing over $\la$, we conclude that
$[Q_3]_p = [V_8]_p - p\sum_{\la=1}^{p-2} \phi(\la) (a_{\la,p}^2-p) - p^2$.
Substituting the expression for $[V_8]_p$ from Theorem~\ref{thm:main-first}
and that for $b_p$ from Proposition~\ref{prop:finish-proof}
completes the proof of Theorem~\ref{thm:q3}.

\subsection{Geometry and Calabi-Yau resolutions}\label{subsec:geom-cy-res}
Having shown that $V_8$ is modular, we now consider the question of whether it
is birationally equivalent to a Calabi-Yau fivefold.  Such problems were
considered by Cynk and Hulek \cite{cynk-hulek}, who gave a sufficient
condition for a double cover of a smooth variety to admit a crepant resolution.
We proceed to state their result and to explain why it does not apply to $V_8$.

\begin{definition}[{\cite[Section~5]{cynk-hulek}}]\label{def:near-pencil}
  Let $\{D_i\}_{i=1}^n$ be a set of smooth divisors on a variety $V$ and
  let $S \subseteq \{1,\dots,n\}$ be a nonempty subset such that
  $\cap_{i \in S} D_i \not \subset D_j$ for all $j \notin S$.
  We say that the intersection
  $\cap_{i \in S} D_i$ is {\em near-pencil} if there is a single element
  $s \in S$ such that $\cap_{i \in S} D_i \ne \cap_{i \in S \setminus \{s\}} D_i$.
\end{definition}

Note in particular that $\cap_{i \in S} D_i$ is automatically near-pencil if
$\#S \le 2$.  For another example, if $V = \P^n$, the $D_i$ are hyperplanes,
and the equation defining $D_1$ involves a variable not mentioned in any
other $D_i$, then every intersection $\cap_{i \in S} D_i$ with $1 \in S$ is
near-pencil.

Cynk and Hulek show:
\begin{proposition}[{\cite[Proposition 5.6]{cynk-hulek}}]\label{prop:cynk-hulek}
  Let $V$ be a smooth variety
  with smooth divisors $D_1, \dots, D_n$ such that the sum of
  the Picard classes of the $D_i$ is divisible by $2$ in $\Pic V$.  For
  $S \subset \{1,\dots,n\}$, let $C_S = \cap_{i \in S} D_i$.  Suppose that, for
  all nonempty $S$ with $C_S \ne C_{S \cup \{i\}}$ for all $i \notin S$, either
  $C_S$ is near-pencil or $\lfloor \#S/2 \rfloor = \codim S - 1$.  Then the
  double cover of $V$ branched along the union of the $D_i$ admits a crepant
  resolution.
\end{proposition}

We refer to the given condition on $S$ or $C_S$ as the
{\em Cynk-Hulek criterion}.  If the condition holds for the
intersection of every subset of the $D_i$ of cardinality greater than $1$,
we will say that the set of $D_i$ satisfies the Cynk-Hulek criterion.

To discuss $V_8$, we do not need to describe the resolution in detail.
%% (however, we will sketch the construction that proves this proposition later,
%% in Method~\ref{crepant-res}).
It suffices to observe that all subsets of
the $12$ hyperplanes satisfy the Cynk-Hulek criterion, with the 
exception of the
intersection of the hyperplanes $x_i + x_{i+1} = 0$, which consists of the
single point $(-1:1:-1:1:-1:1)$ and is not near-pencil (the intersection of
any five of the six hyperplanes is the same).
Combining this with the result of Cynk and Hulek just above, we conclude that
if the singularity of $V_8$ at $(0:-1:1:-1:1:-1:1)$ admits a
crepant resolution, then so does $V_8$, and this resolution $\tilde V_8$
would be a Calabi-Yau fivefold.  However, this does not seem likely.

\begin{conjecture}\label{conj:no-crepant-not-cy}
  The singularity of $V_8$ at $(0:-1:1:-1:1:-1:1)$ does not admit a
  crepant resolution, and hence $V_8$ is not a singular Calabi-Yau fivefold.
\end{conjecture}

A heuristic argument for this statement can be given by using the ideas
of~\cite{c-vs} to calculate the deformation space of a putative Calabi-Yau
resolution $\tilde V_8$.  We expect (but cannot prove) that it satisfies
$h^{4,1}(\tilde V_8) = \sum_{\codim C_i} h^0(K_{C_i})$, where the $C_i$ are the loci
blown up to construct the resolution; on the other hand,
the $\F_p$-point counts suggest that
$h^{4,1}(\tilde V_8) = 1+\sum_{\codim C_i} h^0(K_{C_i})$.

Finally, we use \cite[Theorem 2]{kollar-larsen} to investigate the
resolutions of $Q_3$.

\begin{proposition}\label{prop:age-1}
  The age {\rm (\cite[Definition 1]{kollar-larsen})} of $\iota_3$
  acting on the tangent space of every fixed point is $1$.  
\end{proposition}

\begin{proof}
  Viewed as an automorphism
  of $\P^5$, the fixed locus of $\iota_3$ consists of the linear subspaces
  $x_1 - x_5 = x_2 - x_4 = 0$ and $x_1 + x_5 = x_0 = x_2 + x_4 = x_3 = 0$.
  On the first of these, we may take $x_1, x_2, x_3, x_4, x_5$ as a system of
  local parameters, even on the double cover.  Then $\iota_3$ exchanges the
  tangent vectors in the
  $x_2$ and $x_4$ directions, and likewise $x_1$ and $x_5$, while fixing
  $x_3$; thus its age is $1$ there.
  
  On the second, we have $t = 0$ on the double cover, so $t$ must be taken
  among our local parameters, and we must blow up $x_0 = x_3 = 0$.  Then
  we take $x_0, x_1, x_3, x_4, t$ as our local parameters.  It is clear that
  tangent vectors in the $x_0, x_3, t$ directions are fixed by $\iota_3$.
  As for $x_1$, such a tangent vector is described by the infinitely near
  point $(0:x_0:x_1+\epsilon:-x_4:x_3:x_4:-x_1)$ where $\epsilon^2 = 0$,
  which goes by the involution to $(0:x_0:-x_1:x_4:x_3:-x_4:x_1+\epsilon)$.
  Now $(-x_1:x_1+\epsilon) = (x_1-\epsilon:-x_1)$, so the corresponding
  diagonal entry of the matrix giving the action is $-1$; similarly for a
  tangent vector in the $x_4$ direction.  Thus the action of $\iota_3$ has
  trace $1$ on the $5$-dimensional tangent space, and so the $-1$-eigenspace
  has dimension $2$ and the age is $1$ as before.
\end{proof}

In particular $\iota_3$ satisfies the global Reid-Tai criterion, and so
by \cite[Theorem 2]{kollar-larsen} the quotient would have Kodaira
dimension $0$ if $V_8$ had a Calabi-Yau resolution; such a resolution
would exist if not for the singularity at $(0:-1:1:-1:1:-1:1)$, but
as in Conjecture~\ref{conj:no-crepant-not-cy} we believe that it does not.

Nevertheless, if we blow up only this bad singularity,
the intersections of all subsets of the strict transforms
meet transversely, so the quotient is well-behaved (unfortunately this is
not a crepant blowup).
Because of this and Theorem~\ref{thm:q3}, we conjecture:
\begin{conjecture}\label{conj:rigid-8} $Q_3$ admits a strongly rigid
  Calabi-Yau resolution of singularities for which the representation
  on $H^5$ coincides with that obtained from the newform of weight~$6$ and
  level~$8$ up to semisimplification.
\end{conjecture}

Unfortunately, it seems difficult to study $Q_3$ by the methods of
\cite{cynk-hulek}.  Although $Q_3$ is still a double
cover of a rational variety, the branch locus is now quite complicated,
with some collections of components having nonreduced intersection.
Thus we are unable to prove Conjecture~\ref{conj:rigid-8}; we emphasize,
however, that our proof of modularity is independent of any such
considerations and is unconditional.

%% \begin{remark}\label{rem:q3-hard}
%%   Let us make these statements more precise.  The fixed locus of $\iota_3$
%%   on $\P^5$ consists of two disjoint components, one defined by
%%   $x_0 = x_3 = x_1 + x_5 = x_2 + x_4 = 0$ and one by $x_1 - x_5 = x_2 - x_4 = 0$.
%%   Let these be $C_1, C_2$ respectively.
%%   Blowing them up, we obtain a smooth subvariety of
%%   $\P^5 \times \P^3 \times \P^1$ on which $\iota_3$ is everywhere defined and
%%   whose fixed locus consists of the two exceptional divisors $E_1, E_2$, which
%%   are disjoint.  It is easily checked that the $+1$-eigenspace
%%   for the action on the tangent space has dimension $4$ in each case, so the
%%   quotient is smooth.  In addition, the canonical divisor of the blowup is
%%   $-6H + 3E_1 + E_2$, so the canonical divisor of the quotient is
%%   $-6 \pi_* H + \pi_* E_1$ (the quotient map being ramified along the images
%%   of $E_1$ and $E_2$).  Since $2$ of the branch divisors contain $C_1$ and
%%   none contains $C_2$, our numerical data are still compatible with a crepant
%%   resolution.  However, there are $7$ nonreduced intersection components, which
%%   themselves intersect in various ways.  All of these must be blown up before
%%   we can apply the results of~\cite{cynk-hulek}.
%%   This seems to be well beyond the current limits of computer algebra.
%% \end{remark}

\section{The second example: level $32$}\label{sec:32}
In this section we will consider the fivefold $V_{32}$ defined by the equation
\begin{equation}\label{fivefold-32}
  t^2 = \left(\prod_{i=0}^5 x_i\right)(x_0+x_1)(x_3+x_5)(x_2+x_4+x_5)(x_0+x_2-x_4)(x_1-x_2+x_4)(x_2-x_3+x_4).
\end{equation}
We will show that it realizes the newform of weight~$6$ and level~$32$ that
has complex multiplication by $\Q(i)$ (see Definition~\ref{def:mj-ajp}
for notation):
\begin{theorem}\label{thm:count-32}
  $[V_{32}]_p = \sum_{i-0}^5 p^i - a_{6,p} - pa_{4,p} - 2p^2a_{2,p}$.
\end{theorem}

As with Section~\ref{sec:first}, the assertions of this section
and some that it refers to in Sections
\ref{sec:hypergeom}~and~\ref{sec:geom-points}
are verified numerically in the file {\tt code-level-32.mag} in
\cite{magma-scripts}.  Most of this is routine: for example,
Lemma~\ref{lem:coef-mf} is just a matter of comparing some products of
coefficients and Proposition~\ref{prop:count-k} is a simple application of
Lemma~\ref{lem:count-one-dc}.  However, Proposition~\ref{prop:fibres-rho}
involves counting the points on a toric variety that
is better not embedded in projective space, and
Proposition~\ref{prop:match-fibre-p2p2} requires some care to account for all
of the exceptions in the statement.

\subsection{Proof of modularity}\label{subsec:modularity-32}

We will prove Theorem~\ref{thm:count-32} in a manner suggested by Lemma~\ref{lem:coef-mf}.  Namely, we
will start by writing a fibration on $V_{32}$ whose fibres are quotients of
products of two K3 surfaces.  One of these is always the K3 surface $M$ of
Picard rank $20$ and discriminant $-4$, whose point count over $\F_p$ is
controlled by the $a_{3,p}$ (Lemma~\ref{lem:coef-mf}).
The other varies in a family whose total
space (the $\L$ of Section~\ref{sec:geom-points})
is related to $M \times E_{32}$, where $E_{32}$ is
an elliptic curve of conductor $32$.  Its number of $\F_p$-points is
therefore related to $a_{3,p} a_{2,p}$ and hence to $a_{4,p}$.
Thus $[V_{32}]_p$ involves $a_{4,p} a_{3,p}$ and therefore $a_{6,p}$.

%% \begin{remark}\label{rem:blowing-up-fivefolds-hard}
%%   Proving that a resolution has the expected Hodge numbers is much more
%%   difficult for fivefolds than it is for threefolds.  The problem is that
%%   it could happen that, after some number of blowups, the intersection of
%%   the ramification divisor with the next locus to be blown up is so large
%%   that the branched cover of this locus is not rational.  For threefolds,
%%   blowing up a curve $C$ in the singular locus of the branch divisor
%%   that intersects another
%%   such curve $C'$ removes and adds one point to the intersection of the
%%   ramification divisor with $C'$, which does not cause problems.  However,
%%   blowing up a curve in the singular locus of a fivefold that is contained
%%   in a threefold is a much more delicate matter.
%% \end{remark}

As in Section~\ref{sec:first}, we begin by partitioning
the twelve hyperplanes into the branch locus into two sets of six, each
set intersecting in a line.  In particular, the set of linear forms
$\{x_3,x_5,x_0+x_1,x_3+x_5,x_2+x_4+x_5,x_2-x_3+x_4\}$ spans the space
generated by $x_2+x_4,x_0+x_1,x_3,x_5$, while its complement in the set of
$12$ linear forms defining components of the branch locus spans
$\langle x_0+x_1,x_2+x_4,x_0,x_2 \rangle$.
\begin{definition}\label{def:fib-v32}
  Define a rational map $\pi: V_{32} \dashrightarrow \P^1$ by
  $(x_0+x_1:x_2+x_4)$.  We
  will also view $\pi$ as a map $\P^5 \dashrightarrow \P^1$.
\end{definition} 

As before, the general fibre is a quotient of the product of two K3 surfaces.
As in Definition~\ref{def:fla}, we describe the first of these by
writing the linear form $ax_3 + bx_5 + c(x_0+x_1) + d(x_2+x_4)$ as
$ax + by + (c\la+d)z$.  This gives six linear forms 
$x,y,\la z,x+y,y+z,-x+z$ from which we obtain a K3 surface
$k_\la$ defined by the equation obtained by setting $t^2$ equal to their
product.  Similarly, for the other six we write
$ax_0 + bx_2 + c(x_0+x_1) + d(x_2+x_4)$ as $ax + by + (c\la+d)z$, obtaining
$x,-x+\la z,y,-y+z,x+2y-z,-x-2y+(\la+1)z$ and define a K3 surface
$\ell_\la$ by setting $u^2$ equal to their product.
However, we are only interested in
$(k_\la \times \ell_\la)/\sigma$, where $\sigma$ is the involution that changes
the signs of $t, u$.  This does not change if we replace $\la z$ by $z$
in the definition of $k_\la$ and $y$ by $\la y$ in that of $\ell_\la$.
Thus instead of $k_\la$ and $\ell_\la$, we may use $M, N_\la$, as defined in
Definition~\ref{def:kla-lla-new}.

\begin{definition}\label{def:involutions} Let $E = E_{32}$ be the elliptic curve
  of conductor $32$ whose affine equation is $y^2 = x^3 - x$
  (this was previously referred to as $E_{-1}$, because we wanted to
  emphasize its role in the fibre of $\pi$ at $-1$).
  We use $-$ to denote the automorphism of $K$ defined by
  $(t:x:y:z) \to (-t:x:y:z)$, and also the automorphism of $E_{32}$ taking
  $(x:y)$ to $(x:-y)$.
\end{definition}

\begin{proposition}\label{prop:bir-v32}
  $V_{32}$ is birationally equivalent to
  $(M \times M \times E)/\Sigma$, where
  $\Sigma$ is the group of automorphisms acting as $-1$ on an even number of
  the factors and $+1$ on the remaining ones.
\end{proposition}

\begin{proof} Since $\L$ is the total space of the $N_\la$, our discussion
  above shows that $V_{32}$ is birational to $(M \times \L)/(-,-)$.
  The result follows from Corollary~\ref{cor:totl-is-mtimese-mod-sigma}.
\end{proof}

This will be used in Section~\ref{subsec:rigid-v32}
to construct a rigid Calabi-Yau quotient of $V_{32}$.
It is also instructive to compare to Lemma~\ref{lem:count-q-dc}.
Indeed, $e_{32,+} - e_{32,-} = a_{2,p}$, where $e_{32,\pm}$ are as in
Lemma~\ref{lem:count-q-dc}, and for a suitable 
model of $M$ we have $k_+ - k_- = a_{3,p}$.
%% However, this does not immediately
%% imply the simple formula of Theorem~\ref{thm:count-32}, because the birational
%% equivalence contracts and expands many subvarieties.

First we use Proposition~\ref{prop:bir-v32}
together with Lemma~\ref{lem:count-q-dc} to count the
$\F_p$-points of $V_{32}$.  It follows from Lemma~\ref{lem:count-q-dc} that $(M \times L_\la)/\sigma$ has
$(p^2+p+1)^2 + (P_{M,1}-N_{M,1})(P_{L_\la,1} - N_{L_\la,1})$ points for $\la \in \F_p$.
On the other hand, we may use the birational equivalence of this with the
fibre of $\pi$ at $\la$ to count the points on the fibre.  In the following,
let $y_0,y_1,y_2,z_0,z_1,z_2$ be coordinates on $\P^2 \times \P^2$.

\begin{proposition}\label{prop:match-fibre-p2p2}
  Fix $\la \in \F_p^*$, and let $\mu$ be the rational map from the hyperplane
  in $\P^5$ defined by $x_0+x_1 = \la(x_2+x_4)$ to $\P^2 \times \P^2$ given by
  $((x_3:x_5:x_2+x_4),(x_0:x_2:x_2+x_4))$.  Then $\mu$ induces a rational map
  from the fibre of $\pi$ at $(\la:1)$ to $(M \times N_\la)/\sigma$.
  Further, it induces a bijection between the sets of $\F_p$-points of these
  schemes, with
  the following exceptions:
  \begin{enumerate}
  \item Points of the fibre of $\pi$ with $x_0 = x_1 = x_2 + x_4 = 0$, or
    with $x_3 = x_5 = x_2+x_4 = 0$, do not correspond to any point of
    $(M \times N_\la)/\sigma$.
  \item Points of the fibre of $\pi$ with $x_2 + x_4 = 0$, but
    %$(x_0:x_2) = (x_1:x_4)$, but
    for which $(x_0:x_2), (x_1:x_4), (x_3:x_5)$ are well-defined points of
    $\P^1$, correspond $(p-1)$-to-$1$ to points of
    $(M \times N_\la)/\sigma$ above a point
    with coordinates $((y_3:y_5:0),(y_0:y_2:0))$.
  \item Points of $(M \times N_\la)/\sigma$ with $y_2 = 0, z_2 \ne 0$, or with
    $y_2 \ne 0, z_2 = 0$, do not correspond to any point of the fibre of
    $\pi$.
  \end{enumerate}
\end{proposition}

\begin{proof}
  As discussed above, the map matches the branch loci of the two double
  covers, so there is a rational map from the fibre of $\pi$ to
  $(M \times N_\la)/\sigma$ as described.  In case 1, we would obtain
  a point whose coordinates in one $\P^2$ are $(0:0:0)$.  In case 2,
  if $x_2 + x_4 = 0$, then clearly we obtain the point
  $((x_3:x_5:0),(x_0:x_2:0))$, and this is unchanged by rescaling $(x_0:x_2)$
  by an element of $\F_p^*$.  On the other hand, if $x_0 = x_2 = 0$, then
  we have already dealt with this point in case 1, and similarly for the
  other two pairs.  
  Finally, points with $y_2 = 0, z_2 \ne 0$ or vice versa
  cannot be obtained from $\mu$,
  because $x_2+x_4$ cannot both be $0$ and not be $0$.
  
  In the other direction, we have an inverse rational map from 
  $\P^2 \times \P^2$ to the hyperplane.  On the affine patch
  $y_2 = z_2 = 1$, it is given by
  $((y_0,y_1),(z_0,z_1)) \to (z_0:\la-z_0:z_1:y_0:1-z_1:y_1)$.  Points not
  on this affine patch are accounted for in cases 2 and 3 above.
\end{proof}

\begin{corollary}\label{cor:count-diff}
  The double cover $(M \times N_\la)/\sigma$ of $\P^2 \times \P^2$ has
  $p(p+1)^2+(p-2)\phi(\la)(P_{M,1}-N_{M,1})$ more points than the
  fibre of $\pi$ at $\la$.
\end{corollary}

\begin{proof}
  We consider the three types of exception in
  Proposition~\ref{prop:match-fibre-p2p2}.
  Case 1 describes two disjoint sets of $p+1$
  points in the fibre of $\lambda$,
  so $2p+2$ in total.  In case 2 we contract $(p+1)^2(p-1)$ points to
  $(p+1)^2$.
  
  To understand the third case, note that one of the linear forms defining
  $M$ is the third coordinate, so all points with $x_2 = 0$ give one point
  on $M$ and there are $(p+1)p^2$ missed
  points for which the third coordinate of the point giving $M$ is $0$.
  On the other hand, setting the
  third coordinate to $0$ in the linear forms defining $L_\la$ gives
  $\pm t_0, \la t_1, -t_1, \pm (t_0+2t_1)$.  Thus the product is $0$ for $3$
  points and $\la$ times a nonzero square for the other $p-3$.
  Where the product is $0$, we have $p^2$ points.  
  In the double cover we get one point for every $Z$-point of $M$ (cf.{}
  Proposition~\ref{prop:count-k}) and two points for each $P$-point or
  $N$-point, depending on whether $\la$ is a square.
  This contributes
  $p^2(p+1) + (p-2)(P_{M,1}-N_{M,1})$ points if $\phi(\la) = 1$ and
  $p^2(p+1) + (p-2)(N_{M,1}-P_{M,1})$ points if $\phi(\la) = -1$ to the
  excess of $[(M\times N_\la)/\sigma]_p$ over $[\pi_\la]_p$.
  
  In total, then, the excess is
  $p^2(p+1) + p^2(p+1) + (p-2)\phi(\la)(P_{M,1}-N_{M,1}) - (2p+2) - (p+1)^2(p-1)$.
  This simplifies to the formula asserted.
\end{proof}

These calculations are not valid for $\la = 0, \infty$.  However, it is
easy to see that for both of these the fibre has $p^4+p^3+p^2+p+1$ points.
Finally, the base locus of $\pi$ is defined by $x_0 + x_1 = x_2 + x_4 = 0$.
On this locus the product of linear forms is $0$, so it has
$p^3 + p^2 + p + 1$ points and we must subtract $p$ times this from the
total number of points on the fibres to obtain the correct point count for
$V_{32}$.

We now assemble all of the ingredients: the comparison of the fibres
in Proposition~\ref{prop:match-fibre-p2p2} and its
Corollary~\ref{cor:count-diff}, the remarks on special fibres and
the base locus just above, the count of points on $\L$ from
Proposition~\ref{prop:fibres-rho}, and the relations of coefficients of
modular forms of Lemma~\ref{lem:coef-mf}.
By routine calculation, we obtain
$[V_{32}]_p = \sum_{i-0}^5 p^i - a_{6,p} - pa_{4,p} - 2p^2a_{2,p}$ as
claimed.  This completes the proof of Theorem~\ref{thm:count-32}.

\subsection{Construction of a rigid Calabi-Yau fivefold of level 32}\label{subsec:rigid-v32}
As in Section~\ref{subsec:rigid-f1}, we will use the method of~\cite{burek}
to construct a candidate for a strongly rigid Calabi-Yau fivefold of
level~$32$.
In this case the control given to
us by Proposition~\ref{prop:bir-v32} allows us to construct a strongly rigid
Calabi-Yau fivefold as a quotient of $M \times M \times E$ and to prove
that it is birationally equivalent to the quotient of $V_{32}$ that we
construct here.

\begin{remark}\label{rem:cohom-kke}
  Since $V_{32}$ satisfies the Cynk-Hulek criterion
  (Proposition~\ref{prop:cynk-hulek}), it admits
  a crepant resolution by a Calabi-Yau fivefold.  The shape of the formula
  for the number of points suggests that this resolution has
  $h^{5,0} = h^{4,1} = 1, h^{3,2} = 2$, and $h^{i,j} = 0$ unless
  $i = j$ or $i+j = 5$.
  
  This is explained by the birational description.  Indeed, $H^5$ of the
  resolution arises from $H^2_T(M) \otimes H^2_T(M) \otimes H^1(E)$,
  where $H^2_T$ is the transcendental lattice $H^2(M)/\Pic M$.
  Thus, for example, $H^{3,2}$ of the resolution matches
  $$(H^{2,0}(M) \otimes H^{0,2}(M) \otimes H^{1,0}(E)) \oplus (H^{0,2}(M) \otimes H^{2,0}(M) \otimes H^{1,0}(E))$$
  and has dimension $2$.  This will be
  explained more precisely in the proofs of Theorem~\ref{thm:kke-rigid-quotient}
  and its Corollary~\ref{cor:count-m6}.
\end{remark}

Let $G_{64}$ be the group of projective automorphisms of the configuration
of $12$ hyperplanes used to construct $V_{32}$, let $C_2$ be the cyclic group
of order $2$, and let $Z_G$ be the centre of a group $G$.
The group $G_{64}$ has order $64$ and is isomorphic to
$C_2 \times G_{32}$, where $Z_{G_{32}} \cong C_2^2$ and $G_{32}$ fits into an
exact sequence $1 \to Z_{G_{32}} \to G_{32} \to C_2^3 \to 1$.
% It is (32,27).  Is there a nice description?

We refer to {\tt quotient-level-32.mag} in
\cite{magma-scripts} for verifications of the results of this section.
In particular, we investigate
the automorphisms numerically, verifying that
the point counts of certain quotients of $V_{32}$ are as claimed, thus
justifying our eventual choice of quotient to study more closely.
This is slightly difficult because the branch loci of the quotients are
not always defined by a single equation in our chosen embedding.  Also,
the quotients, though easily embedded in weighted projective space, are
not so pleasant to work with in ordinary projective space.

We study the quotients of $V_{32}$ by elements of order $2$ with
characteristic polynomial $(x-1)^4(x+1)^2$ in $G_{64}$ as in
Section~\ref{subsec:rigid-f1}.
We concentrate on two elements of $G_{64}$, namely
$\alpha_1$, taking $(x_0:x_1:x_2:x_3:x_4:x_5)$ to $(x_1:x_0:x_4 :x_3:x_2:x_5)$,
and %pos2[4] and pos2[2] respectively
$\alpha_2$, defined by $(x_0:x_1:x_2:x_3:x_4:x_5) \to (-x_1:-x_0:x_2:-x_5:x_4:-x_3)$.
(Note that $\alpha_2 \in Z_{G_{64}}$.)  In this case, the Cynk-Hulek criterion
is satisfied, so we know that the differential
\begin{equation}\label{hol-diff-v32}
  D' = \frac{x_5^6}{t} \bigwedge_{i=0}^4 d(x_i/x_5)
\end{equation}
(cf.{} (\ref{hol-diff-f1}), but note that here the argument constitutes
a rigorous proof) pulls back to a generator of $H^{5,0}$ on the
quotient.

Now, $\alpha_1$ gives an even permutation and therefore does
not change the sign of 
$\frac{x_5^5}{\prod_{i=0}^4 x_i} \wedge_{i=0}^4 d(x_i/x_5)$.  In addition, it
fixes $x_5/t$, and so it fixes $D'$.
To verify the invariance for $\alpha_2$, we use the alternative form
$-\frac{x_2^6}{t} \wedge_{\stackrel{i=0}{i \ne 2}}^5 d(x_i/x_2)$ for $D'$.
Negating any variable changes the sign of this, and it is invariant under
even permutations that fix $x_2$.  So it is fixed by $\alpha_2$.

When we consider the quotient $V_{32}/\alpha_1$, we find that the
images of all of the branch divisors are defined by a single polynomial,
and so it is easy to write down the branch function on $\P^5/\alpha_1$
(that is, the function whose square root gives the double cover
$V_{32}/\alpha_1 \to \P^5/\alpha_1$).
On the other hand, for $V_{32}/\alpha_2$, all but two of the branch divisors,
as well as the union of the two that are not, are defined by single
polynomials.  In this case, it is again easy to write down the branch
function.  For both of these, the quotient $\P^5/\alpha_i$ is a hypersurface
of weighted degree $4$ in $\P(1,1,1,1,2,2,2)$ as previously.
There are other involutions such that exactly one branch divisor on the
quotient is not defined by a single polynomial.  We would have to be more
careful in this situation; however, it does not arise in this paper.

Counting $\F_p$-points on the double covers for small odd $p$,
we are led to the following conjecture.
\begin{conjecture}\label{conj:count-mod-a1} For all primes $p>2$ we have
  $[V_{32}/\alpha_1]_p = [V_{32}/\alpha_1 \alpha_2]_p = \sum_{i=0}^5 p^i - a_{6,p} - p^2 a_{2,p}$
  and $[V_{32}/\alpha_2]_p = \sum_{i=0}^5 p^i - a_{6,p} - pa_{4,p}$.
\end{conjecture}

Accordingly we expect that $\alpha_1$ acts as $-1$ on
$H^{4,1}(\tilde V_{32}) \oplus H^{1,4}(\tilde V_{32})$ and has eigenvalues
$1,-1$ on $H^{3,2}(\tilde V_{32}) \oplus H^{2,3}(\tilde V_{32})$, while
$\alpha_2$ acts as $+1$ on $H^{4,1} \oplus H^{3,2}$ and as $-1$ on
$H^{3,2} \oplus H^{2,3}$; this also confirms our observation that both act as
$+1$ on $H^{5,0}$, which also applies to $H^{0,5}$.
Also $\alpha_1 \alpha_2$ should satisfy the same description as $\alpha_1$,
except that the eigenspaces of $\pm 1$ for $H^{3,2}$ and $H^{2,3}$ are reversed.
This is consistent with
the fact that $\alpha_1\alpha_2$ is conjugate to $\alpha_1$ in $G_{64}$.
In particular, the $+1$ eigenspace of $\langle \alpha_1,\alpha_2\rangle$ on
$H^5$ should be neither more nor less than $H^{5,0} \oplus H^{0,5}$, and
we expect that the number of $\F_p$-points on the quotient should be
expressible in terms of powers of $p$, Artin symbols, and $a_{6,p}$ only.

In order to make rigorous computations of cohomology, we will study
$M \times M \times E$ in place of $V_{32}$.  Recall from
Proposition~\ref{prop:bir-v32} that $V_{32}$ is
birationally equivalent to the quotient of $M \times M \times E$
by a Klein four-group.  We will first lift $\alpha_1, \alpha_2$ to
$M \times M \times E$.

\begin{definition}\label{def:kke-vars}
  On $M \times M \times E$, let the coordinates be
  $t_i, x_i, y_i, z_i$ on the $i$th copy of $M$ and $x_3, y_3, z_3$ on $E$.
  Further, let $\omega_{i,j}$ be generators of $H^{j,2-j}$
  of the $i$th copy of $M$ for $j = 0$ or $2$
  and let $\eta_i$ be generators of $H^{i,1-i}(E)$ for $i = 0, 1$.
\end{definition}

\begin{proposition}\label{prop:lift-alpha-kke}
  The automorphisms $\alpha_1, \alpha_2$ lift to automorphisms of
  $M \times M \times E$ of order $4$, given  respectively by
  \begin{align*}
    \bar \alpha_1&: (t_1:x_1:y_1:z_1),(it_2:-z_2:-y_2:-x_2),(-x_3:iy_3:z_3),\\
    \bar \alpha_2&: (it_1:-y_1:-x_1:z_1),(it_2:z_2:y_2:x_2),(x_3:-y_3:z_3).\\
  \end{align*}
\end{proposition}

(There are four choices for each of these lifts, but all are of
order~$4$.  Note also that $\bar \alpha_1^2 = \sigma_1$ and
$\bar \alpha_2^2 = \sigma_2$.)

\begin{proof} 
  The file {\tt kke.mag} in \cite{magma-scripts} defines $M \times M \times E$
  and various quotients inside appropriate toric varieties, constructs the
  rational map $M \times E \dashrightarrow \L$ and hence
  $M \times M \times E \dashrightarrow M \times \L$, defines the
  automorphisms $\alpha_1, \alpha_2$ in coordinates, and verifies that the
  given formulas define automorphisms of $M \times M \times E$ that are
  pullbacks of $\alpha_1, \alpha_2$.
\end{proof}

\begin{definition}\label{def:g4-bar}
  Let $G_4 = \langle \alpha_1, \alpha_2\rangle$; let
  $\bar G_4 = \langle \bar \alpha_1, \bar \alpha_2\rangle$.
\end{definition}

\begin{proposition}\label{prop:kke-cy-quotient}
  The quotient
  $(M \times M \times E)/\bar G_4$
  admits a resolution of singularities $V$ whose canonical divisor is trivial.
\end{proposition}

\begin{proof}
  It is enough to
  show that the Reid--Shepherd-Barron--Tai criterion \cite{reid} is satisfied.
  For all fixed point strata of an element $\beta$ of order
  $n$ in $\bar G_4$
  of codimension greater than $1$, and every $k \in \Z$ with $(k,n) = 1$,
  the action of every power of $\beta$
  on the tangent space at a general point has eigenvalues $\zeta_n^{a_{i,k}}$
  with $0 \le a_i < n$. The criterion requires that $\sum_i a_i \ge n$
  for all $k$ with equality for at least one $k$.
  
  As this is a tedious calculation, we present it only for one element of
  $\bar G_4$, namely $\bar \alpha_1$.  We first consider
  the points $(t_1:x_1:y_1:z_1),(0:1:c:1),(0:0:1)$.
  The tangent space is the direct sum
  of those for $M, M, E$, and $\bar \alpha_1$ preserves this decomposition.
  Clearly the action of $\alpha_1$ on the tangent space to $(t_1:x_1:y_1:z_1)$
  on the first copy of $M$ is trivial.  On the second copy, we have the
  tangent vectors given by $(\epsilon:1:c:1), (0:1:c+\epsilon:1)$, which
  are taken to $(-i\epsilon:1:c:1)$ and $(0:1:c+\epsilon:1)$ respectively,
  indicating an action by the diagonal matrix with entries $-i,1$.  On
  $E$ a tangent vector is given by $(\epsilon:1:0)$, which goes to
  $(i\epsilon:1:0)$, so the action is by $i$.  Thus the $a_i$ are
  $0,0,3,0,1$; since $n = 4$, the result follows in this case.  Similarly,
  considering the points $(t_1:x_1:y_1:z_1),(0:-1:0:1),(0:0:1)$, we find
  eigenvalues $1,1$ on the first copy of $M$ and $i$ on $E$ as before,
  while at $(0:-1:0:1)$ on the second copy of $M$ we find $i, -1$, and
  the $a_i$ are $0,0,1,2,1$.  Replacing $(0:0:1)$ on $E$ by the other fixed
  point $(0:1:0)$ changes nothing, and
  similar methods apply to all other elements of $\bar G_4$.
\end{proof}

We remark that the resolution of $V$ has trivial fundamental group.
%(this is sometimes required in the definition of a Calabi-Yau variety).
To see this, note that
$M \times M \times E$ has fundamental group $\Z^2$, since $\pi_1(M) = 0$.
The action of $\alpha_1$ is given in suitable coordinates by the matrix
$\begin{pmatrix}0&1\\-1&0\end{pmatrix}$, so the fixed subspace is trivial.
Thus $(M \times M \times E)/\bar G_4$ has trivial $\pi_1$, and the
resolution of singularities replaces contractible subvarieties by
rational subvarieties, whose fundamental group is trivial.  Thus by
Seifert-van Kampen $\pi_1(V)$ is likewise trivial.

\begin{theorem}\label{thm:kke-rigid-quotient}
  The resolution $V$ of
  Proposition~\ref{prop:kke-cy-quotient} is a Calabi-Yau variety and
  $H^4(V)$ is generated by fundamental classes of subvarieties.
\end{theorem}

\begin{proof}
  We have shown that $K_V \cong \O_V$, so we need only show that
  $H^{i,0}(V) = 0$ for $1 \le i \le 4$.
  We consider the action of automorphisms on cohomology.  The
  K\"unneth formula determines the cohomology of $M \times M \times E$.
  It is clear that negation acts as $-1$ on $H^{i,2-i}(M)$ for $i = 0, 2$
  and on $H^1(E)$, while acting as $+1$ on $H^{i,i}(M)$ for $0 \le i \le 2$
  and on $H^{i,i}(E)$ for $0 \le i \le 1$.  Hence, for any field $F$ of
  characteristic not equal to $2$, the subring of
  $H^*(M \times M \times E,F)$ fixed by $\sigma_1, \sigma_2$ is generated as
  an $F$-vector space by the fundamental classes of subvarieties and by 
  classes of the form $\omega_{1,j_1} \otimes \omega_{2,j_2} \otimes \eta_k$.
\end{proof}

\begin{corollary}\label{cor:count-m6}
  The Galois representation on $H^5_\etale(V,\Q_\ell)$
  coincides with that attached to $m_6$ (Definition~\ref{def:mj-ajp}).
\end{corollary}

\begin{proof}
  So far we have not used the elements of $\bar G_4$ of order $4$, but
  now we need to.
  First we show that 
  $H^5((M \times M \times E)/\langle \sigma_1, \sigma_2 \rangle$ fixed
  by $\bar G_4$ coincides with
  $H^{5,0} \oplus H^{0,5} (M \times M \times E)$.
  The $E$-component
  of $\bar \alpha_1$ pulls back $dx_3/y_3$ to $d(-x_3)/(iy_3) = i\,dx_3/y_3$,
  so it acts as multiplication by $i$ on $H^{1,0}(E)$.  Continuing in this
  way, we see that $\bar \alpha_1$ acts as $(1,-i,i)$ on the spaces spanned
  by $\omega_{1,1}$, $\omega_{2,1}$, $\eta_1$ respectively, and
  $\bar \alpha_2$ acts as $(i,i,-1)$.  The action on the spaces spanned by
  $\omega_{1,2}$, $\omega_{2,2}$, $\eta_2$ is by the reciprocals of these.

  Thus the $+1$ eigenspace for $\bar G_4$ acting on $H^5$
  is spanned by $\omega_{1,1} \otimes \omega_{2,1} \otimes \eta_1$
  and $\omega_{1,2} \otimes \omega_{2,2} \otimes \eta_2$.  This establishes
  that the dimension of $H^5(V)$ is $2$ and that the Galois representation
  is the component with weights $(5,0)$ and $(0,5)$ in the tensor product
  of those attached to the cusp forms corresponding to $M, M, E$.
  In light of Lemma~\ref{lem:coef-mf}, this is
  the form $m_6$ of weight~$6$ and level~$32$ with complex multiplication
  by $\Q(i)$.
\end{proof}

\begin{corollary}\label{cor:kodaira-dim-v32}
  The Kodaira dimension $\kappa$
  of a resolution of $(M \times M \times E)/\bar G_4$ or of $V_{32}/G_4$ is $0$.
\end{corollary}

\begin{proof}
  For $(M \times M \times E)/G_4$, we have $\kappa \ge 0$,
  since the generator of $H^{5,0}$ survives in the quotient;
  also $\kappa \le 0$, since $\kappa_{M \times M \times E} = 0$.
  It follows that $\kappa_{V_{32}/G_4} = 0$ as well.
\end{proof}

\begin{remark}\label{rem:other-cm-forms} The same method could be used to
  construct rigid Calabi-Yau fivefolds realizing the form of weight~$6$
  and level~$27$: see Example~\ref{ex:level-27} below.
  However, it is clearly insufficient even for
  the form of weight~$6$ with CM by $\Q(\sqrt{-19})$, because there is
  no finite subgroup of $GL_5(\bar \Q)$ with traces in $\Q(\sqrt{-19})$ but
  not in $\Q$.
\end{remark}

\begin{example}\label{ex:level-27}
  Let $K$ be the K3 surface defined by $t^2 = xy(x^4+y^4) + z^6$
  (cf.~\cite[Example 0.4 (11)]{mukai}).  Then $K$ has Picard number $20$
  and admits a non-symplectic automorphism $\omega_K$ of order $3$ given by
  $(t:x:y:z) \to (t:x:y:\zeta_3 z)$.  Let $E$ be the elliptic curve
  $y^2 + y = x^3$ of conductor $27$; it also has an automorphism
  $\omega_E$ of order $3$ given by $(x,y) \to (\zeta_3 x,y)$.  We use
  $-_K, -_E$ for the obvious negation maps; let
  $G \subset \Aut(K \times K \times E)$ be generated by
  $(-_K,-_K,1),(1,-_K,-_E),(\omega_K,\omega_K^{-1},1),(1,\omega_K,\omega_E^{-1})$.
  % the signs are right here: the hol. diff. will be
  % d(x/z) wedge d(y/z) * z^3/t, which pulls back to zeta_3 * same
  Then $(K \times K \times E)/G$ will have a strongly rigid Calabi-Yau
  resolution $C_{27}$.  It is defined over $\Q$, since the group $G$ is
  rational, even though some of its elements are not.
  Since the representation for the modular form
  attached to $K$ is a component of the symmetric square of the form
  attached to $E$,
  that for $C_{27}$ will be a component of $\Sym^5$, which means it must be
  the unique rational newform of weight~$6$ and level~$27$.

  %% Direct attempts to associate rigid Calabi-Yau fivefolds to other forms in
  %% this way are doomed to failure.
  %% As in \cite[Theorem 1.4]{as}, if a K3 surface admits an automorphism
  %% $\alpha$ of finite order such that no power of $\alpha$ other than the
  %% identity is symplectic, then its transcendental lattice has a
  %% fixed-point-free automorphism of the same order.  For a K3 surface of
  %% rank $20$, this means that the automorphism must have order $3$ or $4$
  %% with the transcendental lattice being a multiple of the lattice of rank $2$
  %% and discriminant $3$ or $4$.  The corresponding modular
  %% form of weight~$3$ will always be a quadratic twist of the one of level
  %% $12$ or $16$ respectively.
\end{example}

Since $V_{32}/G_4$ is birationally equivalent to $V$,
we have likewise shown that $V_{32}$ admits a rigid Calabi-Yau resolution.
In fact, considering that these are related by blowing up and down minimal
rational varieties defined at worst over $\Q(i)$, we conclude:

\begin{proposition}\label{prop:count-quo-v32}
  There is a formula of the form $\sum_{i=0}^5 c_i(p) p^i - a_{6,p}$ for the
  number of $\F_p$-rational points of $V_{32}/G_4$, where $c_0 = c_5 = 1$ and
  the other $c_i$ are expressions of the form $c + c' \phi(-1)$ with
  $c, c' \in \frac12 \Z$.
\end{proposition}

%% \begin{remark}\label{rem:cohom-v32}
%%   Although $H^5(K \times K \times E)$ is much larger, the classes coming
%%   from elements of $\Pic K$ do not survive in the quotient by the involutions
%%   $\sigma_1, \sigma_2$ (for notation see Corollary~\ref{prop:bir-v32}).
%%   This is because they are fixed by these involutions,
%%   while $H^1(E)$ and the transcendental part of one $H^2(K)$ are negated by
%%   $\sigma_1$.  Thus $H^3(\L)$ has dimension $4$ and is negated by
%%   $\sigma_2$ (which acts as $-1$ on $H^1(E)$ and $+1$ on the same $H^2(K)$),
%%   while the other $H^2(K)$ is negated by $\sigma_2$ as well.  So we find
%%   $8$ for the dimension of $H^5$ of a resolution of $V_{32}$, and
%%   $H^{5,0}$ arises from $H^{1,0}(E) \otimes H^{2,0}(K) \otimes H^{2,0}(K)$,
%%   etc.  In particular $H^{3,2} \oplus H^{2,3}$ comes from
%%   $H^1(E) \otimes (H^{2,0}(K) \oplus H^{0,2}(K))^{\otimes 2}$, which is isomorphic
%%   to the sum of two copies of $H^1(E)$ twisted by $2$, which is why we see
%%   $2p^2a_{2,p}$ in Theorem~\ref{count-32}.  Similarly $H^{4,1} \oplus H^{1,4}$
%%   corresponds to the Hecke character $\chi^4 \bar \chi$ and its conjugate.
%%   Since $\chi \bar \chi$ takes $p$ to $p$ when $p$ is a prime congruent to
%%   $1 \bmod 4$, we obtain a $pa_{4,p}$ term.
%% \end{remark}

\begin{theorem}\label{thm:rigid-32}
  $[V_{32}/G_4]_p = \sum_{i=0}^5 p^i - a_{6,p}$ for all $p>2$.
\end{theorem}

\begin{proof}
  To determine the $c_i$ of Proposition~\ref{prop:count-quo-v32},
  we consider $V_{32}/G_4$, and in particular the ring
  of invariants of $G_4$.  By calculation in Magma we find that
  the primary invariants have degrees $1,1,2,2,2,2$, while the secondary
  invariants have degrees $0,2,2,4$, with the secondary invariant of degree
  $4$ being the product of the two of degree $2$.  So the quotient is
  embedded in weighted projective space $\P(1,1,2,2,2,2,2,2)$, or (composing
  with the embedding by $\O(2)$) in $\P^8$.
  The branch locus has $5$ orbits of size $2$ and one of size $1$
  under the group; the images of the components in orbits of size $2$ are
  defined by a single polynomial, as is the union of the two of size $1$.
  As before we are able to compute the number of $\F_p$-points of
  $[V_{32}/G_4]$ for small $p$, finding it to be
  $\sum_{i=0}^5 p^i - a_{6,p}$ for $p < 25$.  There are $8$ odd primes less than
  $25$ and this gives us $8$ independent conditions on the coefficients, so
  by basic linear algebra we conclude that this formula holds for all $p$.
\end{proof}

\bibliography{art}

\begin{thebibliography}{10}
\bibitem{ahlgren}
  S. Ahlgren.  {\em The points of a certain fivefold over finite fields
    and the twelfth power of the eta function.}  Finite Fields and their
  Applications 8 (2002), 18--33.
  \url{https://doi.org/10.1006/ffta.2001.0315}.

%% \bibitem{aop}
%%   S. Ahlgren, K. Ono, D. Penniston.  {\em Zeta functions of an infinite family
%%   of K3 surfaces.}  Amer. J. Math 124 (2002), 353--368.
%%   \url{https://doi.org/10.1353/ajm.2002.0007}.
  
%% \bibitem{as} M. Artebani, A. Sarti.  {\em Non-symplectic automorphisms
%%   of order $3$ on K3 surfaces.}  Math. Ann. 342 (2008), 903--921.
%%   \url{https://doi.org/10.1007/s00208-008-0260-1}.
  
\bibitem{bcdt}
  C. Breuil, B. Conrad, F. Diamond, R. Taylor.  {\em On the modularity of
    elliptic curves over $\Q$: Wild $3$-adic exercises.}  J. Amer. Math. Soc.
  14(4) (2001), 843--939.
  \url{https://doi.org/10.1090/S0894-0347-01-00370-8}.
  
\bibitem{magma}
  W. Bosma, J. Cannon, and C. Playoust.  {\em The Magma algebra system.  I.
    The user language}, J. Symbolic Comput. 24 (1997), 235--265.
  
\bibitem{burek}
  D. Burek.  {\em Rigid realizations of modular forms in Calabi–Yau threefolds.}
  J. Pure Appl. Algebra 223 (2019), 547--552.
  \url{https://doi.org/10.1016/j.jpaa.2018.04.005}.
  
\bibitem{cynk-hulek}
  S. Cynk, K. Hulek.  {\em Higher-dimensional
    modular Calabi-Yau manifolds.}  Canad. Math. Bull. 50 (4), 2007, 486--503.
  \url{https://doi.org/10.4153/CMB-2007-049-9}.
  
%% \bibitem{cynk-meyer}
%%   S. Cynk, C. Meyer.  {\em Geometry and arithmetic of certain double octic
%%     Calabi-Yau manifolds.}  Canad. Math. Bull. 48 (2), 2005, 180--194.

\bibitem{c-sz}
  S. Cynk, T. Szemberg.  {\em Double covers of $\P^3$ and Calabi-Yau
    varieties.}  Banach Center Publications 44 (1), 1998, 93--101.
  %\url{arXiv:math/9902057}.
  \url{https://doi.org/10.4064/-44-1-93-101}.

\bibitem{c-vs}
  S. Cynk, D. van Straten, {\em Infinitesimal deformations of double covers
    of smooth algebraic varieties.}  Math. Nachr. 279 (7), 716--726.
  \url{https://doi.org/10.1002/mana.200310388}.

\bibitem{deligne}
  P. Deligne.  {\em Formes modulaires et r\'epresentations $\ell$-adiques.}
  S\'eminaire Bourbaki 355 (1968-69), 139--172.
  \url{http://www.numdam.org/item?id=SB_1968-1969__11__139_0}.
  
\bibitem{dieu-man}
  L. Dieulefait, J. Manoharmayum.  {\em Modularity of rigid Calabi-Yau threefolds
    over $\Q$.}  In {\em Calabi-Yau varieties and Mirror Symmetry}
  (N. Yui, J. D. Lewis, eds.), 159--166.  American Mathematical Society, 2003.
  \url{http://dx.doi.org/10.1090/fic/038}.
  
\bibitem{elkies-schutt}
  N. Elkies, M. Sch\" utt.  {\em Modular forms and K3 surfaces.}
  Adv. Math. 240 (2013), 106--131.
  \url{https://doi.org/10.1016/j.aim.2013.03.008}.
  
\bibitem{fop}
  S. Frechette, K. Ono, M. Papanikolas.  {\em Gaussian hypergeometric functions
    and traces of Hecke operators.}  Int. Math. Res. Not. 60 (2004), 3233--3262.
  \url{https://doi.org/10.1155/S1073792804132522}.
  
  %% \bibitem{vg-top} B. van Geemen, J. Top.  {\em An isogeny of K3 surfaces.}
  %%   Bull. London Math. Soc. 38 (2006), 209--223.
  
\bibitem{gouvea-yui}
  F. Gouv\^ea, N. Yui.  {\em Rigid Calabi-Yau threefolds over $\Q$ are modular.}
  Exp. Math. 29 (1) (2011), 142--149.
  \url{https://doi.org/10.1016/j.exmath.2010.09.001}.
  
  %% \bibitem{greene}
  %%   J. Greene.  {\em Hypergeometric functions over finite fields.}  Trans. AMS
  %%   301 (1) (1987), 77--101.
  
\bibitem{h} R. Hartshorne.  {\em Algebraic geometry.}  GTM 52.  Springer-Verlag, 1977.
  
%% \bibitem{ingalls-logan}
%%   C. Ingalls, A. Logan.  {\em On the Cynk-Hulek criterion for crepant resolutions of double covers.}  \url{https://arxiv.org/abs/math/2006.14981},
%%   submitted.
  
%% \bibitem{ireland-rosen}
%%   K. Ireland, M. Rosen.  {\em A classical introduction to modern number theory}
%%   (second edition).  GTM 84.  Springer-Verlag, 1990.

%% \bibitem{kk}
%%   J. Keum, S. Kondo, {\em The automorphism groups of Kummer surfaces associated with the product of two elliptic curves}, Trans. Amer. Math. Soc. 353, no. 4 (2001), 1469--1487.  \url{https://www.jstor.org/stable/221865}.

\bibitem{kollar-larsen}
  J. Koll\'ar, M. Larsen, {\em Quotients of Calabi-Yau varieties.}
  In {\em Algebra, Arithmetic, and Geometry} (Y. Tschinkel, Y. Zarhin, eds.).
  Volume II, 179--211.  Progress in Mathematics 270.  Birkh\" auser, 2010.
  %\url{https://arxiv.org/abs/math/0701466}.
  \url{https://doi.org/10.1007/978-0-8176-4747-6_6}.
  
%% \bibitem{ks}
%%   M. Kuwata, T. Shioda, {\em Elliptic parameters and defining equations
%%     for elliptic fibrations on a Kummer surface.}  In {\em Algebraic geometry
%%     in East Asia---Hanoi 2005} (K. Konno, V. Nguyen-Khac, eds.).
%%   Adv. Stud. Pure Math 50 (2008), 177--215.
%%   \url{https://doi.org/10.2969/aspm/05010177}.
  
\bibitem{llt}
  W.-C. W. Li, L. Long, F.-T. Tu, {\em A Whipple $_7F_6$ formula revisited},
  \url{https://arxiv.org/pdf/2103.08858.pdf}.

\bibitem{lmfdb}
  The LMFDB Collaboration, {\em The $L$-functions and modular forms database},
  \url{http://www.lmfdb.org}.
  
\bibitem{magma-scripts}
  A. Logan.  Magma scripts.  \url{https://arxiv.org/e-print/1811.02739}.
  
\bibitem{meyer}
  C. Meyer.  {\em A dictionary of modular threefolds.}  PhD thesis, 
  Johannes Guten\-berg-Uni\-ver\-sit\" at in Mainz, 2005.
  \url{https://d-nb.info/975958038/34}.

\bibitem{mukai} S. Mukai, {\em Finite groups of automorphisms of K3 surfaces
  and the Mathieu group}, Inv. Math. 94 (1988), 183--221.
  \url{https://doi.org/10.1007/BF01394352}.
  
\bibitem{p-r} K. Paranjape, D. Ramakrishnan, {\em Modular forms and Calabi-Yau varieties.} In {\em Arithmetic and geometry} (L. Dieulefait, D. R. Heath-Brown, G. Faltings, Yu. I. Manin, B. Z. Moroz, and J.-P. Wintenberger, eds.), 351--372.  LMS Lecture Note Series 420.  London Mathematical Society, 2015.
  %\url{arXiv:1404.1154}.
  \url{https://doi.org/10.1017/CBO9781316106877.019}.

\bibitem{reid} M. Reid, {\em Canonical 3-folds.} In {\em Journ\'ees de
  G\'eometrie Alg\'ebrique d’Angers} (A. Beauville, ed.), 273--310.
  Sijthoff \& Nordhoff, 1980.
  \url{https://doi.org/10.2969/aspm/00110131}.
  
\bibitem{ribet} K. Ribet, {\em Galois representations attached to eigenforms with Nebentypus.}  In {\em Modular functions of one variable, V} (J.-P. Serre and D. Zagier, eds.), 18--52.  Lecture Notes in Mathematics 601.  Springer, 1977.
  \url{https://doi.org/10.1007/BFb0063943}.
  
\bibitem{roberts} D. Roberts, {\em Newforms with rational coefficients},
  Ramanujan J. 46 (3) (2018), 835--862.
  \url{https://doi.org/10.1007/s11139-017-9914-5}.
  
\bibitem{schutt} M. Sch\" utt, {\em CM newforms with rational coefficients.}
  Ramanujan J. 19 (2009), 187--205.
  \url{https://doi.org/10.1007/s11139-008-9147-8}.

%% \bibitem{ss} M. Sch\" utt, T. Shioda, {\em Elliptic surfaces}, \url{arXiv:0907.0298}.
\end{thebibliography}

\end{document}